\documentclass[12pt,reqno]{amsart}

\usepackage[T1]{fontenc}
\usepackage{lmodern}
\usepackage{amsmath,amssymb,mathtools}
\usepackage[a4paper,left=44.69mm,textwidth=5in,top=40.7mm,
textheight=584pt,headheight=8pt,headsep=14pt,footskip=12pt]{geometry}
\usepackage{booktabs,array,needspace,graphicx,enumitem,microtype,etoolbox}
\usepackage[colorlinks=true,linkcolor=blue,citecolor=blue,urlcolor=blue]{hyperref}
\usepackage[nameinlink,noabbrev]{cleveref}
\AtBeginEnvironment{thebibliography}{\raggedright}

\makeatletter
\patchcmd{\@settitle}{\bfseries}{\normalfont}{}{}
\patchcmd{\section}{\scshape}{\upshape}{}{}
\patchcmd{\subsection}{\bfseries}{\normalfont}{}{}
\makeatother

\newtheoremstyle{template}{8pt}{8pt}{\normalfont}{}{\bfseries}{.}{0.5em}{}
\theoremstyle{template}
\newtheorem{maintheorem}{Theorem}

\newtheorem{theorem}{Theorem}
\newtheorem{proposition}[theorem]{Proposition}
\newtheorem{lemma}[theorem]{Lemma}
\newtheorem{corollary}[theorem]{Corollary}

\newtheorem{remark}[theorem]{Remark}
\newtheorem{conjecture}[theorem]{Conjecture}
\crefname{maintheorem}{Theorem}{Theorems}
\Crefname{maintheorem}{Theorem}{Theorems}

\newcommand{\B}{\mathbb B}
\newcommand{\C}{\mathbb C}

\newcommand{\wt}{\widetilde}

\newcommand{\ip}[2]{\langle #1,#2\rangle}

\title[Fixed-radius averages and Toeplitz compactness]
{Compact Toeplitz Operators and the Geometry of Fixed-Radius Averages}
\author[Pan Ma]{Pan Ma}
\address{School of Mathematics and Statistics, HNP-LAMA, Central South
	University, Changsha 410083, China}
\email{pan.ma@csu.edu.cn}
\author[Zhaoxi Zhang]{Zhaoxi Zhang}
\address{School of Mathematics and Statistics, HNP-LAMA, Central South
	University, Changsha 410083, China}
\email{csuzhangzhaoxi@163.com}
\date{}
\hypersetup{pdftitle={Fixed-Radius Averages and Toeplitz Compactness: Strongly Pseudoconvex Domains and the Polydisc},pdfauthor={}}
\begin{document}
\begin{abstract}
	
We characterize compact Toeplitz operators with bounded complex-valued
symbols on the Bergman space of a bounded strongly pseudoconvex domain
with $C^2$ boundary. For any fixed $R>0$, the operator $T_f$ is compact
if and only if the Euclidean-volume averages of $f$ over Bergman
geodesic balls of radius $R$ tend to zero at the boundary. This answers
Zhu's fixed-radius question affirmatively on the unit ball.
On the polydisc $\mathbb D^n$, $n\ge2$, averages over products of
Bergman discs of fixed radii characterize compactness, whereas
averages over genuine Bergman geodesic balls need not.
For every sufficiently small fixed radius, we construct a bounded
real-valued symbol, smooth in the domain, whose geodesic-ball averages
tend to zero along every sequence approaching the boundary, while
its Berezin transform converges to a positive value along one such
sequence. Boundary scaling reduces the strongly pseudoconvex case
to a uniqueness problem for fixed-radius ball averages on complex
hyperbolic space. We prove the required uniqueness theorem under a
sharp decay condition expressed in terms of a Busemann height function.

\end{abstract}
\maketitle
\markboth{Fixed-radius averaging and Toeplitz compactness}
{Fixed-radius averaging and Toeplitz compactness}

\section{Introduction}
The theorem of Axler--Zheng~\cite{AxlerZheng} characterizes
compactness of finite sums of finite products of Toeplitz operators
with bounded symbols on the Bergman space of the unit disk by the
boundary vanishing of the Berezin transform.  Engli\v{s}~\cite{Englis1999} obtained related results on bounded
symmetric domains. Su\'arez~\cite{Suarez} and
Mitkovski--Su\'arez--Wick~\cite{MSW} established compactness criteria for the norm-closed Toeplitz algebras on the unit ball and its
standard weighted Bergman spaces and
Mitkovski--Wick~\cite{MW2014} developed the corresponding theory
on the polydisc. On smoothly bounded strongly pseudoconvex domains,
Wang--Xia~\cite[Theorem~1.2]{WX} proved that an operator
in the Toeplitz algebra generated by bounded symbols is
compact whenever its Berezin transform vanishes at the boundary.

For a single Toeplitz operator, the Berezin transform averages
its symbol against the squared moduli of normalized Bergman
kernels.
This
suggests replacing these kernel averages by Euclidean-volume averages
over Bergman geodesic balls of one fixed radius.
In his talk \emph{Open Problems on Bergman Spaces},
Zhu~\cite{ZhuOpen2020} proposed a conjecture whose fixed-radius
part can be stated as follows.

\begin{conjecture}[Zhu, fixed-radius part]
	Let $f\in L^\infty(\mathbb B_n)$ and fix $R>0$.
	The Toeplitz operator $T_f$ is compact on $A^2(\mathbb B_n)$
	if and only if
	\[
	\frac{1}{V(D_{\mathbb B_n}(z,R))}
	\int_{D_{\mathbb B_n}(z,R)} f(w)\,dV(w)
	\longrightarrow 0
	\qquad (|z|\to1),
	\]
	where $D_{\mathbb B_n}(z,R)$ denotes the Bergman geodesic
	ball centered at $z$ with radius $R$,
	$dV$ is Euclidean volume and $V(E)=\int_E dV$.
\end{conjecture}

 Cao--He--Zhang~\cite[Theorem~1.1]{CHZ} proved this criterion on the
unweighted unit disk. Hagger--Liu--Taskinen--Virtanen~\cite{HLTV} studied sufficient
conditions for boundedness and compactness in terms of local
averages for more general locally integrable symbols.

We prove the fixed-radius criterion on every bounded strongly
pseudoconvex domain with $C^2$ boundary, thereby answering the
question affirmatively on the unit ball. We also show that the
criterion depends on the averaging sets. On the polydisc
$\mathbb D^n$, $n\ge2$, products of fixed Bergman discs characterize
compactness. For every sufficiently small fixed radius, however,
there is a bounded smooth real-valued symbol whose geodesic-ball
averages vanish along the entire boundary while its Toeplitz operator
is not compact. Geometric inclusions and volume comparisons alone
cannot transfer the criterion between the two families, because
the averages retain cancellation.

We now introduce the notation needed to state the main results.

Let $\Omega\Subset\mathbb C^n$ be a connected domain and let $dV$ denote
Euclidean volume. Write $P_\Omega$ for the orthogonal projection onto
$A^2(\Omega)=L^2(\Omega,dV)\cap\mathcal O(\Omega)$, and put
\[
T_f=P_\Omega M_f\big|_{A^2(\Omega)},\qquad
\widetilde f(z)=\langle T_fk_z,k_z\rangle,
\qquad f\in L^\infty(\Omega),
\]
Here $M_f$ denotes multiplication by $f$, $K_\Omega$ is the Bergman
kernel, and $k_z=K_\Omega(\cdot,z)/K_\Omega(z,z)^{1/2}$.
Symbols are complex-valued unless otherwise specified. For a bounded operator $S$, write
$\|S\|_{\mathrm{ess}}=\inf\{\|S-K\|:K\text{ is compact}\}$.
For definiteness we normalize the
Bergman metric by
\begin{equation}\label{eq:metric-convention}
	g_\Omega=\frac1{n+1}\operatorname{Re}
	\bigl(\partial\bar\partial\log K_\Omega(z,z)\bigr).
\end{equation}
Thus the squared length of $v\in\mathbb C^n$ at $z$ is
\[
g_\Omega(z)(v,v)=\frac1{n+1}
\sum_{j,k=1}^n
\frac{\partial^2\log K_\Omega(z,z)}{\partial z_j\partial\bar z_k}
v_j\overline{v_k}.
\]
The unit ball then has holomorphic sectional curvature $-4$.
Any other fixed positive normalization only rescales the radii below.
Let $D_\Omega(z,R)$ be the geodesic ball of radius $R$ and define
\begin{equation*}
	\mathcal A_R^\Omega f(z)
	=\frac{1}{V(D_\Omega(z,R))}
	\int_{D_\Omega(z,R)}f(w)\,dV(w).
\end{equation*}
These averages use the Euclidean measure defining $A^2(\Omega)$,
rather than the Riemannian volume of $g_\Omega$.

Our first main result identifies compactness with boundary vanishing
of either the Berezin transform or the fixed-radius averages.

\begin{maintheorem}
\label{thm:psc-fixed-radius}
Let $n\ge1$ and let the entire boundary of $\Omega\Subset\mathbb C^n$
be $C^2$ and strongly pseudoconvex. For $f\in L^\infty(\Omega)$ and
each fixed $R>0$, the following conditions are equivalent:
\begin{enumerate}[label=\textup{(\roman*)}]
\item $T_f$ is compact on $A^2(\Omega)$;
\item $\widetilde f(z)\to0$ as $z\to\partial\Omega$;
\item $\mathcal A_R^\Omega f(z)\to0$ as $z\to\partial\Omega$.
\end{enumerate}
Thus vanishing for one positive radius is equivalent to vanishing for
every positive radius.
\end{maintheorem}

All boundary limits are taken along arbitrary sequences approaching
the topological boundary. In particular, on the polydisc we require
a limit along the entire boundary, not just the distinguished boundary.

In complex dimension one, the Levi condition on the complex tangent
space is vacuous. The theorem therefore includes bounded $C^2$
planar domains with several boundary components.

On the polydisc, we distinguish products of one-dimensional metric
balls from balls for the Riemannian product metric.
Let $D_{\mathbb D}(a_j,r_j)$ denote a one-dimensional Bergman disc, put
\[
Q_{\mathbf r}(a)=\prod_{j=1}^nD_{\mathbb D}(a_j,r_j),\qquad
\mathcal Q_{\mathbf r}f(a)
=\frac{1}{V(Q_{\mathbf r}(a))}
\int_{Q_{\mathbf r}(a)}f\,dV,
\]
and continue to use $\mathcal A_R^{\mathbb D^n}$ for a genuine
geodesic ball of the product Bergman metric.

The next theorem gives an exact criterion for product averages and
a counterexample for geodesic-ball averages.

\begin{maintheorem}
\label{thm:intro-polydisc}
Let $n\ge2$.
\begin{enumerate}[label=\textup{(\roman*)}]
\item For $f\in L^\infty(\mathbb D^n)$ and every fixed vector
$\mathbf r=(r_1,\ldots,r_n)$ with $r_j>0$, the operator $T_f$ is
compact if and only if $\mathcal Q_{\mathbf r}f(a)\to0$ as
$a\to\partial\mathbb D^n$.
\item There is $R_n>0$ such that, for each $0<R<R_n$, there are
$t_R>0$ and a real-valued symbol $f_R\in C^\infty(\mathbb D^n)$,
$|f_R|\le1$, satisfying
\[
 \begin{gathered}
\mathcal A_R^{\mathbb D^n}f_R(a)\longrightarrow0
 \quad(a\to\partial\mathbb D^n),
 \\
 \|T_{f_R}\|_{\mathrm{ess}}\ge
 (1+t_R^2)\left(\frac{\pi t_R}{\sinh(\pi t_R)}\right)^2>0.
 \end{gathered}
\]
The averages use Euclidean volume and genuine Bergman geodesic
balls, and the limit is along the entire topological boundary.
\end{enumerate}
\end{maintheorem}

The proof of Theorem~\ref{thm:psc-fixed-radius} uses holomorphic
boundary scaling. Along any sequence approaching the boundary,
one can pass to a subsequence for which the rescaled symbols
converge locally weak-$*$ to a bounded function $F$ on a Siegel
domain. The limiting Bergman metric is complex hyperbolic, and the
rescaled balls converge to balls for this metric. The kernel and
metric convergence are obtained from boundary localization and
scaling as in Kim--Yu~\cite[Theorem~4 and Section~3]{KimYu1996};
the precise statements used here are collected in
Proposition~\ref{prop:scv-scaled-geometry}.

The central step is to show that vanishing fixed-radius averages
force every such limit $F$ to be zero. In the Siegel realization,
Euclidean and invariant volume satisfy
$dV_{\mathrm{Euc}}=c_n\rho^{n+1}dV_{\mathrm{CH}}$, where
$\rho$ is the Busemann height at the distinguished point at infinity.
Thus $H=\rho^{n+1}F$ has zero invariant-volume integral over every
ball of the prescribed radius and satisfies
$|H|\le\|f\|_\infty\rho^{n+1}$.
Theorem~\ref{sp:one-radius} proves uniqueness under the more general
bound $|H|\le C\rho^m$ for sufficiently small positive $\rho$,
provided $m>n/2$. No bound is imposed where $\rho$ is large, and
nonzero smooth solutions exist for $m\le n/2$.
Since $n+1>n/2$, the theorem gives $F=0$ almost everywhere.
A localized estimate of the Toeplitz form then yields compactness.
This argument applies at $C^2$ boundary regularity.

On the polydisc, products of Bergman discs allow successive
applications of the one-dimensional uniqueness argument.
For genuine geodesic balls, the product metric couples the
coordinate distances, allowing cancellation between different
coordinates. We construct an oscillatory model with vanishing
ball averages and a nonvanishing Berezin transform. A localization
argument preserves the vanishing of the averages along the entire
topological boundary while retaining a strictly positive lower
bound for the essential norm of the associated Toeplitz operator.

The complex hyperbolic uniqueness theorem is proved using the
Darboux identity and a radial Neumann expansion, together with
local boundary uniqueness for Laplace eigenfunctions.
Fixed-radius mean equations and spectral expansions on rank-one
symmetric spaces were studied by
Volchkov~\cite{Volchkov2001,Volchkov2002}; see also
Platonov~\cite{Platonov2000}.
The eigenfunction uniqueness statement used here is deduced from
the local boundary-value theory of
Kashiwara--\={O}shima~\cite{KO1977}, as developed in
Hansen--Hilgert--Parthasarathy~\cite[Section~6 and Appendix~A]{HHP2019}
and \={O}shima~\cite{Oshima1983}.
Its local formulation allows us to impose the height bound only
where the Busemann height is small, without a growth restriction
at the distinguished boundary point.
This hypothesis differs from the global weighted-growth condition
in Volchkov~\cite[Theorem~1]{Volchkov2001} and the boundary-region
and horizontal-decay assumptions in
Ochakovskaya~\cite[Theorem~3]{Ochakovskaya2013}.

\subsection*{Organization}

Section~\ref{sp:section} proves the complex hyperbolic uniqueness
statement. Section~\ref{sec:strongpsc-fixed-radius} establishes the
strongly pseudoconvex theorem, and Section~\ref{pd:scope} treats product
averages and the geodesic-ball counterexample on the polydisc.
Section~\ref{sec:weighted-ball} proves the weighted ball criterion and
gives examples of compactness caused by oscillation. The final section
discusses shrinking radii and boundary approaches. The appendices
supply a Fourier--Hermite proof of model rigidity and an example
clarifying the boundary regularity assumption.

\section{A model uniqueness theorem}
\label{sp:section}

In this section we prove the model uniqueness theorem used in the
boundary scaling argument. We use the complex hyperbolic metric on
$X=\mathbb{CH}^n$ with holomorphic sectional curvature $-4$.
The standard geometry and Darboux identities used below are described in
Helgason~\cite{Helgason}; the Siegel normalization agrees with
Lu--Yang~\cite{LuYang}. Its Laplace--Beltrami operator $\Delta_X$ has the sign convention
$\operatorname{spec}(-\Delta_X)=[n^2,\infty)$.
For $\zeta\in\partial_\infty X$, write $P_\zeta$ for the positive
Busemann height, normalized by $P_\zeta(o)=1$ at a chosen base point.
In the unit-ball realization this is
\[
 P_\zeta(z)=\frac{1-|z|^2}{|1-\langle z,\zeta\rangle|^2}.
\]
The normalization of the metric gives
\begin{equation}
 \Delta_X(P_\zeta^s)=4s(s-n)P_\zeta^s,
 \qquad |\nabla\log P_\zeta|=2.
 \label{sp:height-identities}
\end{equation}
For complex $s$, the power means
$P_\zeta^s=\exp(s\log P_\zeta)$, using the real logarithm of the
positive height. In particular,
\begin{equation}
 e^{-2R}P_\zeta(x)\le P_\zeta(y)\le e^{2R}P_\zeta(x)
 \quad\text{if }d(x,y)\le R.
 \label{sp:height-comparison}
\end{equation}

The bound below is imposed only towards the finite part of the boundary
in the Siegel realization. No bound is required where the height is large.

\begin{theorem}[One-radius uniqueness under a height bound]
\label{sp:one-radius}
Fix $R>0$, $\zeta\in\partial_\infty X$, and $m>n/2$.
Suppose $H\in L^\infty_{\mathrm{loc}}(X)$ satisfies
\begin{equation}
 \int_{B(x,R)}H(y)\,dV(y)=0\qquad(x\in X)
 \label{sp:zero-ball-means}
\end{equation}
and, for some $C>0$ and $p_0>0$,
\begin{equation}
 |H(x)|\le C P_\zeta(x)^m
 \quad\text{for almost every }x\text{ with }0<P_\zeta(x)\le p_0.
 \label{sp:one-ended-bound}
\end{equation}
Then $H=0$ almost everywhere.
For every $m\le n/2$, there is a nonzero smooth function satisfying
\eqref{sp:zero-ball-means} and \eqref{sp:one-ended-bound}, even with
$C=p_0=1$. Thus the strict inequality $m>n/2$ is sharp.
\end{theorem}

The proof uses Proposition~\ref{sp:eigenfunction-rigidity} and a
radial Neumann expansion. The height-power examples establish the
sharpness of the exponent.

\subsection{A consequence of classical boundary uniqueness}

We first deduce the eigenfunction statement from local
regular-singular boundary-value theory.

\begin{proposition}[Consequence of classical boundary uniqueness]
\label{sp:eigenfunction-rigidity}
Let $\nu>n^2$ and $m>n/2$.
Suppose $u\in C^\infty(X)$ satisfies
\[
 \Delta_Xu=-\nu u
\]
and, for some $C,p_0>0$,
\[
 |u(x)|\le CP_\zeta(x)^m
 \qquad\bigl(0<P_\zeta(x)\le p_0\bigr).
\]
Then $u=0$ on $X$.
\end{proposition}

\begin{proof}
Use the unit-ball realization, with the origin as base point, and put
\[
 \sigma=\sqrt{\nu-n^2}>0,
 \qquad y=e^{-r},\qquad r=d_X(o,x).
\]
The parameter $y$ is the radial boundary parameter; it is distinct
from the Siegel height $\rho$, a positive multiple of $P_\zeta$,
used in the alternative proof.
For $b\in\partial\mathbb B^n$ we have
\[
 \begin{gathered}
x=t(y)b,\qquad t(y)=\frac{1-y^2}{1+y^2},
 \\
 P_\zeta(t(y)b)
 =\frac{4y^2}
 {(1+y^2)^2\lvert1-t(y)\langle b,\zeta\rangle\rvert^2}.
 \end{gathered}
\]
Choose a nonempty boundary patch
$V\Subset\partial\mathbb B^n\setminus\{\zeta\}$.
The last formula gives $P_\zeta\asymp_V y^2$ on a sufficiently small
collar over $V$.  Shrinking the collar if necessary, the hypothesis
therefore implies
\[
 |u(b,y)|\le C_Vy^{2m},
 \qquad
 \sup_{b\in V}\lvert y^{-n\pm i\sigma}u(b,y)\rvert
       \le C_Vy^{2m-n}\longrightarrow0.
\]
Both normalized functions extend continuously by zero at $y=0$.
In particular, all distributions needed below are extendible on this
local collar.

The characteristic exponents of the eigenvalue equation in the
parameter $y$ are $s_\pm=n\pm i\sigma$.
The local boundary-value construction for a regular-singular equation
identifies a continuous normalized limit with its boundary value;
see \cite[Lemma~6.1 and Appendix~A, Lemma~A.1]{HHP2019}.
It follows that both local boundary values vanish on $V$.
Here we use the construction on germs over $V$, rather than a global
boundary-value theorem on the whole boundary.
The two characteristic exponents are nonresonant.  For $n\ge2$, the
root normalization in that reference has parameters
$\varrho=n$ and $\lambda=i\sigma$.
For $n=1$, its indivisible-root convention instead uses
$Y=y^2$, $\varrho=1/2$, and $\lambda=i\sigma/2$;
the normalized functions $Y^{-\varrho\pm\lambda}u$ are exactly
$y^{-1\pm i\sigma}u$.  Thus the same conclusion holds in dimension one.

The local two-boundary-value uniqueness theorem now makes $u$ zero
on an interior neighborhood of any smaller nonempty boundary patch.
For this implication, see the microlocal identity and the
Holmgren--Kashiwara argument in \cite[Remark~6.3, equation~(28)]{HHP2019},
which attributes this uniqueness statement to
\cite[Proposition~4.9]{KO1977}; see also the boundary-value construction
in \cite{Oshima1983}.
All these operations take place over $V$, so no bound or extendibility
at $\zeta$ is required.  Finally, interior unique continuation for
$(\Delta_X+\nu)u=0$ on the connected manifold $X$ gives $u=0$ everywhere.
\end{proof}

\subsection{The positive radial Neumann spectrum}

Let $\omega_{2n-1}$ be the area of the unit sphere in
$T_oX$.
In geodesic polar coordinates, write
\begin{equation*}
 \begin{gathered}
dV=J_n(r)\,dr\,d\omega,
 \qquad
 J_n(r)=(\sinh r)^{2n-1}\cosh r,
 \\
 \Delta_{\mathrm{rad}}=
       \frac1{J_n(r)}\frac d{dr}\left(J_n(r)\frac d{dr}\right),
 \end{gathered}
\end{equation*}
where $d\omega$ is the usual angular area measure.
In particular, $J_n(r)=r^{2n-1}(1+O(r^2))$ at zero.
Multiplying $dV$ and $J_n$ by a common positive constant has no effect
on any argument below.

The radial subspace reduces the Neumann Laplacian on the compact
geodesic ball $\overline{B(o,R)}$.
It therefore has an orthonormal eigenbasis
$q_k$ in $L^2((0,R),J_n(r)\,dr)$, with
\begin{equation}
 -(J_nq_k')'=\nu_kJ_nq_k,
 \qquad q_k'(R)=0,
 \qquad 0=\nu_0<\nu_1\le\nu_2\le\cdots,
 \label{sp:neumann-system}
\end{equation}
and $\nu_k\to\infty$.
These are smooth radial functions at the center, so $q_k'(r)=O(r)$.
The function $q_0$ is constant.

\begin{lemma}
\label{sp:neumann-gap}
Every positive radial Neumann eigenvalue in
\eqref{sp:neumann-system} satisfies $\nu_k>n^2$.
\end{lemma}

\begin{proof}
Suppose $0<\nu\le n^2$ and choose
\[
 s=\frac{n+\sqrt{n^2-\nu}}2\in\mathbb R.
\]
By \eqref{sp:height-identities}, $P_\zeta^s$ is a strictly positive
solution of $\Delta_Xv=-\nu v$.
Average it over the isotropy group of $o$, using probability Haar measure.
The resulting radial function $v_\nu$ is strictly positive, satisfies
$v_\nu(0)=1$, and solves
\[
 (J_nv_\nu')'+\nu J_nv_\nu=0.
\]
Hence
\[
 J_n(R)v_\nu'(R)
       =-\nu\int_0^R J_n(r)v_\nu(r)\,dr<0.
\]
It cannot satisfy a Neumann condition at $R$.
Every nonzero smooth radial solution of this equation is a multiple
of the regular solution with value $1$ at zero, by uniqueness of the
regular radial initial-value problem.
Thus no positive radial Neumann eigenvalue lies in $(0,n^2]$.
\end{proof}

\subsection{Derivation of the model statement}

\begin{proof}[Proof of Theorem~\ref{sp:one-radius}]
First suppose $H$ is smooth. Write
\[
 M_rH(x)=\frac1{\omega_{2n-1}}
          \int_{S_xX}H(\exp_x(r\omega))\,d\omega.
\]
On the rank-one symmetric space $X$, the Darboux identity is
\begin{equation}
 \Delta_xM_rH(x)=
       \frac1{J_n(r)}\partial_r\bigl(J_n(r)\partial_rM_rH(x)\bigr).
 \label{sp:darboux}
\end{equation}
Moreover,
\begin{equation}
 M_rH(x)=H(x)+O_x(r^2),\qquad
 \partial_rM_rH(x)=O_x(r)
 \quad(r\downarrow0).
 \label{sp:center-behavior}
\end{equation}
The zero-ball-mean assumption reads
\begin{equation}
 \int_0^R M_rH(x)J_n(r)\,dr=0.
 \label{sp:radial-zero-mean}
\end{equation}
Apply $\Delta_x$ and use \eqref{sp:darboux}.
The lower endpoint vanishes by \eqref{sp:center-behavior}, giving
\begin{equation*}
 \partial_rM_rH(x)\big|_{r=R}=0.
\end{equation*}

For the basis in \eqref{sp:neumann-system}, set
\[
 c_k(x)=\int_0^R M_rH(x)q_k(r)J_n(r)\,dr.
\]
Equation \eqref{sp:radial-zero-mean} gives $c_0=0$.
For $k\ge1$, integration by parts in \eqref{sp:darboux} gives
\begin{equation}
 \Delta_Xc_k=-\nu_kc_k.
 \label{sp:coefficient-eigenfunction}
\end{equation}
The boundary term is
$[J_n((\partial_rM_rH)q_k-M_rHq_k')]_0^R$.
It vanishes at $R$ by the two Neumann conditions, and at zero because
it is $O_x(r^{2n})$.

If $P_\zeta(x)\le p_0e^{-2R}$, every point of $B(x,R)$ lies in the
region where \eqref{sp:one-ended-bound} applies.
Equation \eqref{sp:height-comparison} then implies
\[
 |c_k(x)|\le C_{k,R}P_\zeta(x)^m
 \qquad\bigl(P_\zeta(x)\le p_0e^{-2R}\bigr).
\]
By Lemma~\ref{sp:neumann-gap}, $\nu_k>n^2$.
Proposition~\ref{sp:eigenfunction-rigidity} applied to
\eqref{sp:coefficient-eigenfunction} therefore gives $c_k=0$.
Completeness of the radial Neumann basis shows that
$M_rH(x)=0$ in $L^2((0,R),J_n(r)\,dr)$ for each $x$.
Continuity in $r$, followed by $r\downarrow0$, gives $H(x)=0$.

For general $H\in L^\infty_{\mathrm{loc}}(X)$, identify $X$ with the
solvable Iwasawa group $NA$, acting simply transitively by isometries
and fixing $\zeta$.
Choose a compactly supported smooth approximate identity
$\eta_\varepsilon$ on this group, with left Haar measure $dg$, and set
\[
 H_\varepsilon(x)=
       \int_{NA}\eta_\varepsilon(g)H(g^{-1}x)\,dg.
\]
This is the usual smooth group convolution; it is valid also though
$NA$ is not unimodular.
The isometric action and Fubini's theorem preserve
\eqref{sp:zero-ball-means}.
For each $g\in NA$, the ratio
$P_\zeta(g^{-1}x)/P_\zeta(x)$ is a positive constant independent of $x$.
It is bounded above and below as $g$ varies on the compact support
of $\eta_\varepsilon$.
Thus $H_\varepsilon$ satisfies a bound of the form
\eqref{sp:one-ended-bound}, possibly with a smaller positive height
cutoff $p_{0,\varepsilon}$ and a different constant.
The smooth argument gives $H_\varepsilon=0$.
As $\varepsilon\downarrow0$, group convolution converges in $L^1$
on compact subsets, proving $H=0$ almost everywhere.

Finally, choose any $k\ge1$ and put
\[
 \tau_k=\frac12\sqrt{\nu_k-n^2}>0,
 \qquad s_k=\frac n2+i\tau_k,
 \qquad H_k=P_\zeta^{s_k}.
\]
By \eqref{sp:height-identities}, $\Delta_XH_k=-\nu_kH_k$.
The spherical mean of any smooth eigenfunction with this eigenvalue
solves the regular radial equation with initial value equal to the
eigenfunction's value at its center.
Consequently
\[
 M_rH_k(x)=H_k(x)\,\frac{q_k(r)}{q_k(0)}.
\]
Here $q_k(0)\ne0$ by uniqueness of the regular radial solution.
Integrating \eqref{sp:neumann-system} gives
\[
 \int_0^R q_k(r)J_n(r)\,dr
       =-\frac{J_n(R)q_k'(R)}{\nu_k}=0.
\]
Hence $H_k$ has zero radius-$R$ ball means.
It is nonzero and satisfies
\[
 |H_k(x)|=P_\zeta(x)^{n/2}
       \le P_\zeta(x)^m
 \qquad(0<P_\zeta(x)\le1)
\]
for every $m\le n/2$.
This proves the asserted sharpness.
If a real-valued example is desired, use
$P_\zeta^{n/2}\cos(\tau_k\log P_\zeta)$ instead.
\end{proof}

\begin{remark}
\label{sp:scope-remark}
For bounded-symbol boundary profiles in the unweighted Bergman
space, the height exponent is $n+1$, strictly above the threshold
$n/2$. The eigenfunction argument uses local boundary data away
from $\zeta$, so no bound at large height is needed.
\end{remark}

\section{Strongly pseudoconvex domains}
\label{sec:strongpsc-fixed-radius}

Throughout this section, $\Omega\Subset\mathbb C^n$ has $C^2$ strongly
pseudoconvex boundary. We prove Theorem~\ref{thm:psc-fixed-radius}.
Let
$\delta(z)=\operatorname{dist}(z,\partial\Omega)$ in a boundary collar.
The metric is normalized as in \eqref{eq:metric-convention}.

\begin{proposition}[Uniform holomorphic boundary scalings]
\label{prop:scv-scaled-geometry}
Put
\[
 \begin{gathered}
\mathcal H=\{w\in\mathbb C^n:
           \rho(w)=-2\operatorname{Re}w_1-|w'|^2>0\},\\
 o=(-1/2,0),\quad
 g_{\rm CH}=\operatorname{Re}\partial\bar\partial(-\log\rho).
 \end{gathered}
\]
For every boundary sequence $a_j\to\partial\Omega$ there are
holomorphic polynomial automorphisms $G_j$ of $\mathbb C^n$, with
constant complex Jacobian, such that $G_jo=a_j$ and, writing
$\delta_j=\delta(a_j)$ and $\Omega_j=G_j^{-1}(\Omega)$,
\begin{align}
 |J_{\mathbb C}G_j|^2&\asymp\delta_j^{n+1},
       \label{eq:geom-jacobian}\\
 \delta_j^{-1}\delta\circ G_j&\longrightarrow\rho
             \quad\hbox{in }C^2_{\rm loc}(\mathcal H),
       \label{eq:geom-height}\\
 K_{\Omega_j}&\longrightarrow K_{\mathcal H}
             \quad\hbox{in }C^\infty_{\rm loc}
                    (\mathcal H\times\mathcal H),
       \label{eq:geom-kernel}\\
 G_j^*g_\Omega&\longrightarrow g_{\rm CH}
             \quad\hbox{in }C^\infty_{\rm loc}(\mathcal H).
       \label{eq:geom-metric}
\end{align}
Every model compact set is eventually contained in $\Omega_j$.
All assertions are uniform in moving boundary centers. Moreover,
\begin{equation}\label{eq:geom-log-height}
 |d\log\delta|_{g_\Omega}\le C
\end{equation}
in a smaller collar; in particular the metric is complete and proper.
\end{proposition}

\begin{proof}
For a compact $C^2$ boundary, the inner signed distance is $C^2$
on a common tubular collar. Indeed, the nearest-point projection
$\pi$ is $C^1$, and $\nabla\delta=-\nu\circ\pi$, where $\nu$ is the
outward $C^1$ unit normal. Thus $\nabla\delta$ is $C^1$.

Let $q_j$ be the nearest boundary point to $a_j$. A complex affine
map centered at $q_j$, followed by a triangular quadratic change of
the complex normal coordinate involving horizontal variables only,
gives a holomorphic polynomial automorphism $T_j$ such that
\[
 \begin{gathered}
T_j(q_j)=0,\quad T_j(a_j)=(-\delta_j/2,0),\\
 \delta(T_j^{-1}\zeta)
 =-2\operatorname{Re}\zeta_1-|\zeta'|^2
  +Q_j(\zeta)+E_j(\zeta).
 \end{gathered}
\]
Here $Q_j$ is a real quadratic polynomial, each of whose monomials
contains at least one real component of the complex normal variable
$\zeta_1$. The function $E_j$ is the second-order Taylor remainder.
There is a common modulus $\omega(t)\to0$ as $t\downarrow0$ such that
\begin{equation}\label{eq:geom-peano}
 |D^\beta E_j(\zeta)|\le
 C|\zeta|^{2-|\beta|}\omega(C|\zeta|),\qquad |\beta|\le2,
\end{equation}
uniformly in $j$ on a common neighborhood of zero.
All derivatives here are real derivatives. This is the uniform
$C^2$ Taylor formula.
The horizontal linear map normalizes the positive Levi form.
The normal affine map sends the inward unit normal to $(-1/2,0)$.
The quadratic correction vanishes on the normal axis, so the formula
for $T_j(a_j)$ is exact. The coefficients and inverses of these
normalizations have uniform bounds as $q_j$ ranges over the compact
boundary: the normal frames and Levi inverse square roots may be
chosen continuously on finitely many boundary charts.
Uniform continuity of the original second derivatives, and these
polynomial coordinate changes, give the common modulus in
\eqref{eq:geom-peano}. Define
\[
 D_j(w)=(\delta_jw_1,\sqrt{\delta_j}\,w'),\qquad
 G_j=T_j^{-1}\circ D_j.
\]
The Jacobians of the triangular changes are constant, proving
\eqref{eq:geom-jacobian}. We verify the derivative count for
\eqref{eq:geom-height}. Fix a model compact set $L\Subset\mathcal H$.
On $L$ one has $|D_jw|\le C_L\sqrt{\delta_j}$.
A derivative with $a$ normal and $b$ horizontal real derivatives,
$a+b\le2$, applied to $\delta_j^{-1}E_j(D_jw)$, contributes the
factor $\delta_j^{-1+a+b/2}$. By \eqref{eq:geom-peano} its absolute
value is at most
\[
 C_L\delta_j^{a/2}\omega(C_L\sqrt{\delta_j}).
\]
Each normal-containing monomial in $\delta_j^{-1}Q_j(D_jw)$ has a
factor $\sqrt{\delta_j}$ or $\delta_j$, also after two derivatives
in $w$. Consequently
\begin{equation}\label{eq:geom-C2-estimate}
 \left\|\delta_j^{-1}\delta\circ G_j-\rho\right\|_{C^2(L)}
 \le C_L\left[\sqrt{\delta_j}+
                    \omega(C_L\sqrt{\delta_j})\right]\longrightarrow0.
\end{equation}
This proves \eqref{eq:geom-height}, and eventual inclusion of every
compact subset of $\mathcal H$ follows as well. In dimension one
there are no horizontal variables; the same estimates remain valid.

We derive the kernel limit from boundary localization. A strongly pseudoconvex $C^2$ boundary point is a
local holomorphic peak point. Kim--Yu
\cite[Theorem~4, pp.~145--148]{KimYu1996} localize the minimum
integrals on a bounded pseudoconvex domain at such a point. The
zeroth integral is the reciprocal diagonal Bergman kernel, so for a
fixed small neighborhood $U$ of the point,
\[
 \frac{K_\Omega(z,z)}{K_{\Omega\cap U}(z,z)}\longrightarrow1
\]
as $z$ tends to the boundary there, without an approach restriction.
Their Section~3.1 proposition and Section~3.2 moving-center scaling
argument give convergence of the scaled local kernels to
$K_{\mathcal H}$ on compact subsets, including all interior
derivatives. This is the bounded-domain squeeze and normal-family
argument after a fractional change from the Siegel model to the
ball.

We use auxiliary local convexifying coordinates to obtain the
local kernel limit.
In the $T_j$ coordinates, put $s=\operatorname{Im}\zeta_1$.
On the real tangent space $\operatorname{Re}\zeta_1=0$, the
quadratic part of the inner defining function is
\[
 -|\zeta'|^2+b_js^2+2s\ell_j(\zeta'),
\]
where the real coefficients of $b_j$ and of the real-linear form
$\ell_j$ are uniformly bounded. A common sufficiently large
$A>0$ makes the local holomorphic map
\[
 \Psi_j(\zeta)=(\zeta_1+A\zeta_1^2,\zeta')
\]
convexifying: its inverse changes the $s^2$ coefficient to
$b_j-2A$, so the real second fundamental form with the outer
orientation is uniformly positive. The uniform $C^2$ remainder
then gives locally strongly convex domains in a common smaller
neighborhood. These maps are biholomorphic on a common smaller neighborhood. The scaled transition is
\[
 \begin{gathered}
H_j=D_j^{-1}\circ\Psi_j\circ D_j,
 \qquad H_j(w)=(w_1+A\delta_jw_1^2,w'),
 \\
 J_{\mathbb C}H_j=1+2A\delta_jw_1.
 \end{gathered}
\]
Thus $H_j\to\mathrm{id}$ and $J_{\mathbb C}H_j\to1$ in
$C^\infty$ on model compact sets, and
$H_j(o)-o=(A\delta_j/4,0)$.
The auxiliary local kernel limit therefore pulls back to the
original constant-Jacobian maps $G_j$. Only those original maps
are used for the Euclidean-volume averages.

Diagonal localization transfers the local conclusion to the global
kernel. More explicitly, the kernels
\[
 \mathcal K_j(w,v)=J_{\mathbb C}G_j(w)
       \overline{J_{\mathbb C}G_j(v)}
       K_\Omega(G_jw,G_jv)=K_{\Omega_j}(w,v)
\]
are locally bounded on product compact sets: use diagonal domain
monotonicity and kernel Cauchy--Schwarz. Their diagonal converges to
$K_{\mathcal H}(w,w)$ by localization and the local scaling result.
Polarization determines every normal-family limit from its diagonal;
Cauchy estimates give \eqref{eq:geom-kernel} at every derivative
order. These arguments are uniform in moving boundary centers:
otherwise a violating sequence has a boundary-convergent subsequence
to which the same localization and scaling apply.

Since $K_{\mathcal H}(w,w)=c_n\rho^{-n-1}$, taking mixed second
derivatives of its logarithm proves \eqref{eq:geom-metric}.
The holomorphic Jacobian term has zero mixed derivative.
At a normalized center, write
$\delta\circ G_j=\delta_j\rho_j$, with $\rho_j(o)=1$.
Equation~\eqref{eq:geom-height} bounds $d\rho_j(o)$, and
\eqref{eq:geom-metric} bounds the inverse pulled-back metric there.
This proves \eqref{eq:geom-log-height} at every collar center.
A finite-length curve cannot approach the boundary, since
$\log\delta\to-\infty$ there. Standard completeness and the
Hopf--Rinow theorem imply properness; the remaining interior
compact sets cause no difficulty.
\end{proof}

\begin{remark}[Boundary derivatives and interior derivatives]
\label{rem:two-boundary-derivatives}
The estimate \eqref{eq:geom-C2-estimate} uses only continuity of the
second derivatives of a defining function. The arbitrarily high
\emph{interior} derivative convergence in
\eqref{eq:geom-kernel}--\eqref{eq:geom-metric} follows instead from
holomorphy and Cauchy estimates. It does not require corresponding
boundary derivatives. The ball confinement, normalized kernel mass,
and Toeplitz form estimates below use only locally uniform metric
convergence, a first-order height bound, and constant Jacobians.
Thus they impose no additional boundary regularity.
\end{remark}

\subsection{Boundary profiles and normalized kernel mass}

We first specify the geometric inputs supplied by
Proposition~\ref{prop:scv-scaled-geometry}.  Put
\[
 \begin{gathered}
\mathcal H=\{w=(w_1,w'):-2\operatorname{Re}w_1-|w'|^2>0\},
 \\ \rho=-2\operatorname{Re}w_1-|w'|^2,\qquad
 b_0=(-1/2,0).
 \end{gathered}
\]
Given $a_j\to\partial\Omega$ and $r_j=\operatorname{dist}
(a_j,\partial\Omega)$, pass to a subsequence and choose the
normal-axis-preserving holomorphic polynomial normalizations
$G_j=T_j^{-1}\circ D_j$, where
$D_j(w_1,w')=(r_jw_1,\sqrt{r_j}\,w')$.
They are global biholomorphisms of $\mathbb C^n$,
$G_jb_0=a_j$, and their complex Jacobians $J_j$ are constant, with
$|J_j|^2\asymp r_j^{n+1}$.
For $\Omega_j=G_j^{-1}(\Omega)$, every compact subset of $\mathcal H$
is eventually contained in $\Omega_j$, and
\begin{equation*}
 K_{\Omega_j}\longrightarrow K_{\mathcal H}
       \quad\hbox{in }C^\infty_{\rm loc}
              (\mathcal H\times\mathcal H),\qquad
 G_j^*g_\Omega\longrightarrow g_{\rm CH}
       \quad\hbox{in }C^\infty_{\rm loc}.
\end{equation*}
Only fixed compact subsets of these holomorphic charts are used below.

A \emph{boundary profile} of $f$ is any weak-$*$ limit on compact
subsets of $\mathcal H$ of $f\circ G_j$.  Such limits exist after
subsequencing and satisfy $\|F\|_\infty\leq\|f\|_\infty$.
Equivalently one can extend the functions by zero off $\Omega_j$
in $\mathcal H$ and use weak-$*$ compactness in
$L^\infty(\mathcal H)$.  We say that $f$ has zero boundary profiles
if every such limit is zero.

\begin{lemma}[Tightness and the Berezin criterion for profiles]
\label{lem:psc-profile-berezin}
Normalized kernels tend weakly to zero when their centers tend to
$\partial\Omega$.  If $f\circ G_j\to F$ weak-$*$ locally, then, for
each fixed $b\in\mathcal H$,
\begin{equation}\label{eq:psc-berezin-profile}
 \widetilde f(G_jb)\longrightarrow
       \int_{\mathcal H}F(w)|k_b^{\mathcal H}(w)|^2\,dV_{\rm Euc}(w).
\end{equation}
Consequently, $\widetilde f\to0$ at the boundary if and only if
$f$ has zero boundary profiles.
\end{lemma}

\begin{proof}
Let $k_{j,b}$ be the normalized kernel of $\Omega_j$ centered at
$b$.  By the kernel transformation law, the Berezin transform on
the left of \eqref{eq:psc-berezin-profile} is exactly
\[
 \int_{\Omega_j}(f\circ G_j)(w)|k_{j,b}(w)|^2\,dV_{\rm Euc}(w).
\]
Local kernel convergence implies $k_{j,b}\to k_b^{\mathcal H}$
uniformly on each model compact set.  Both kernels have total
$L^2$ mass one.  For any compact $L\Subset\mathcal H$,
\[
 \int_{\Omega_j\setminus L}|k_{j,b}|^2
       \longrightarrow
       1-\int_L|k_b^{\mathcal H}|^2.
\]
Exhausting $\mathcal H$ by compact sets makes the right-hand side
arbitrarily small.  Thus no kernel mass is lost outside the model
region, including the part of $\Omega_j$ outside $\mathcal H$.
Weak-$*$ convergence on $L$ and the uniform bound on $f$ now prove
\eqref{eq:psc-berezin-profile}.

The same argument proves weak convergence to zero in the original
space.  Given a compact $L$ capturing all but an arbitrarily small
kernel mass, $G_j(L)$ has Euclidean volume tending to zero because
$|J_j|^2\asymp r_j^{n+1}$.  For $v\in A^2(\Omega)$, absolute
continuity makes $\int_{G_j(L)}|v|^2\to0$; Cauchy--Schwarz on this
set and its complement proves
$\langle v,k_{G_jb}\rangle\to0$.
Taking $b=b_0$ treats every boundary sequence.

We also recall why the model Berezin transform is injective on
bounded symbols.  If the right-hand side of
\eqref{eq:psc-berezin-profile} vanishes for every $b$, consider
the holomorphic--antiholomorphic kernel
\[
 A_F(z,v)=\int_{\mathcal H}
 K_{\mathcal H}(z,w)F(w)K_{\mathcal H}(w,v)\,dV_{\rm Euc}(w).
\]
It satisfies
$A_F(b,b)=K_{\mathcal H}(b,b)\widetilde F(b)=0$.
Polarization therefore gives $A_F\equiv0$, and hence $T_F=0$.
Under a Cayley biholomorphism from the ball to $\mathcal H$ the
unitary kernel-Jacobian map carries this to $T_{F\circ C}=0$ on
the ball.  Testing against all holomorphic monomials shows
\[
 \int_{\mathbb B_n}(F\circ C)(z)z^\alpha\bar z^\beta\,dV(z)=0
 \quad\hbox{for all }\alpha,\beta.
\]
Stone--Weierstrass, followed by uniqueness of finite measures, gives
$F=0$ almost everywhere.

For each fixed $b\in\mathcal H$, the height convergence gives
\[
 \delta(G_jb)=r_j\bigl(\rho(b)+o(1)\bigr)\longrightarrow0.
\]
Thus the boundary limit applies at every center $G_jb$ in the
same scaling sequence, without a new normalization.
If $\widetilde f\to0$, apply
\eqref{eq:psc-berezin-profile} at every fixed $b$ to obtain this
injectivity conclusion for every profile.  Conversely, if the
Berezin transform failed to tend to zero, choose offending
centers, extract a profile, and apply the same formula at $b_0$.
Zero profiles give a contradiction.
\end{proof}

\subsection{Convergence of metric balls and one-radius totality}

\begin{lemma}[Fixed-radius means and profiles]
\label{lem:psc-profile-means}
For each fixed $R>0$, $\mathcal A_R^\Omega f\to0$ at the boundary
if and only if $f$ has zero boundary profiles.
\end{lemma}

\begin{proof}
Write $g_j=G_j^*g_\Omega$.  For fixed $b\in\mathcal H$ and
$0<\eta<R$, local metric convergence gives, eventually,
\begin{equation}\label{eq:psc-ball-sandwich}
 B_{\rm CH}(b,R-\eta)
       \subset B_{g_j}(b,R)
       \subset B_{\rm CH}(b,R+\eta).
\end{equation}
Here the last inclusion requires a brief confinement argument.
The closed model ball of radius $R+\eta$ is compact, since
complex hyperbolic space is complete and proper.  A path from $b$
to its complement must first cross its boundary.  On that compact
ball the metrics are uniformly close, so the part of the path up
to its first exit has $g_j$ length at least
$(1-o(1))(R+\eta)>R$.  Paths outside the normalization window
therefore cannot produce a shortcut.  The first inclusion follows
by using model minimizing geodesics and the upper comparison of
the two metrics.  This proves \eqref{eq:psc-ball-sandwich}.

Model spheres have Euclidean measure zero.  It follows that the
ball indicators converge in $L^1$, with all the balls confined
to a common model compact set.  Since $G_j$ has constant
Jacobian, it cancels from the normalized Euclidean mean.
For any locally weak-$*$ convergent profile sequence,
\begin{equation}\label{eq:psc-mean-profile}
 \mathcal A_R^\Omega f(G_jb)\longrightarrow
 \frac{\displaystyle\int_{B_{\rm CH}(b,R)}
                   F\,dV_{\rm Euc}}
      {\displaystyle V_{\rm Euc}(B_{\rm CH}(b,R))}.
\end{equation}
Thus zero profiles imply vanishing means, by the usual
contradiction subsequence argument.

For the converse, vanishing means make every numerator in
\eqref{eq:psc-mean-profile} zero, for all model centers $b$.
The model volume relation is
\[
 dV_{\rm Euc}=c_n\rho^{\,n+1}dV_{\rm CH}.
\]
Consequently $H=\rho^{n+1}F$ has zero invariant ball integrals
of radius $R$ and satisfies
$|H|\leq\|f\|_\infty\rho^{n+1}$.
Apply the model uniqueness statement
(Theorem~\ref{sp:one-radius}) at the ideal endpoint whose
Busemann height is $\rho$.  Its required exponent is strictly
above $n/2$, whereas here $n+1>n/2$.  It gives $H=0$, hence $F=0$.
\end{proof}

\subsection{Vanishing boundary profiles imply compactness}

We now show that vanishing boundary profiles imply compactness.  Let $\Omega\Subset
\mathbb C^n$ be a strongly pseudoconvex domain with $C^2$ boundary,
let $g$ be a fixed normalization of its Bergman metric, and let $d$
be the associated distance.  Write $dV$ for Euclidean volume and
$d\mu_g$ for Riemannian volume.  Take $\delta(q)=\operatorname{dist}(q,\partial\Omega)$ throughout
$\Omega$; differentiability of $\delta$ is used only in a boundary
collar.

Here are the scaling data used in this subsection.  The usual
holomorphic boundary scalings, with a fixed model point $o$, give
embeddings $G_q$ centered at $q$ and defined on every fixed compact
subset of the model domain when $q$ is sufficiently close to the
boundary.  On such compact sets assume
\begin{equation}\label{eq:pz-scaling-data}
 \begin{gathered}
 G_q(o)=q,\qquad G_q^*dV=j_q\,dV_{\rm mod},\qquad
 j_q\asymp\delta(q)^{n+1},\\
 G_q^*g\longrightarrow g_{\rm mod}
       \quad\hbox{locally uniformly}.
 \end{gathered}
\end{equation}
where $j_q$ is a positive constant in the profile variables and
$g_{\rm mod}$ is a smooth positive metric.  All assertions are
uniform in the boundary centers.  In the application the model is
the Siegel realization of complex hyperbolic space.  Its exact
curvature is not needed in the present subsection.  We use also
the usual defining-function control
$|d\log\delta|_g\leq C$ in a boundary collar.  It was proved directly in Proposition~\ref{prop:scv-scaled-geometry}.
In particular $d$ is proper: a curve tending to the boundary has
infinite length, and the remaining closed subsets are compact.
No rate of convergence in \eqref{eq:pz-scaling-data} is used.

\begin{lemma}[A single-symbol criterion from boundary profiles]
\label{lem:profile-zero-compactness}
Assume the scaling data just stated.  Let $f\in L^\infty(\Omega)$
be complex-valued.  Suppose that every boundary scaling sequence
has only the zero weak-$*$ profile: for every $q_j\to\partial\Omega$,
every subsequential weak-$*$ limit of $f\circ G_{q_j}$ on compact
subsets of the model domain is zero.  Then the Toeplitz operator
$T_f=P_\Omega M_f$ is compact on $A^2(\Omega)$.
\end{lemma}

\begin{proof}
We first give uniform local geometric estimates.  Choose a small
Euclidean coordinate ball about $o$ whose slightly larger concentric
ball is relatively compact in the model domain.  On this larger
ball, \eqref{eq:pz-scaling-data} bounds the profile metric above and
below by fixed positive multiples of the Euclidean metric.  The
same is true for the profile metric volume density.  A curve from
$o$ to the boundary of the coordinate ball has a fixed positive
minimum length.  Consequently a sufficiently small metric ball
about $q$ cannot leave the image of that coordinate ball: otherwise
its initial segment up to the first exit would already have this
minimum length.  Straight coordinate segments give the converse
inclusion for smaller Euclidean balls.  Comparing the resulting
Euclidean balls and their volume densities proves that there are
$t_*>0$ and $c,C>0$ such that
\begin{equation}\label{eq:pz-local-volumes}
 \begin{gathered}
c t^{2n}\leq\mu_g(B_d(q,t))\leq C t^{2n},\\
 cj_qt^{2n}\leq V(B_d(q,t))\leq Cj_qt^{2n}
 \quad(0<t\leq t_*).
 \end{gathered}
\end{equation}
For centers in a fixed compact interior set use ordinary
holomorphic coordinate balls and take $j_q=1$ there; enlarging the
constants makes all these estimates global.  The choices for
interior centers serve only the estimates and are unrelated to
the boundary-profile hypothesis.

Fix a sufficiently small constant $a>0$, independent of all later
parameters, and put $Q(q)=B_d(q,a)$.  The coordinate balls in the
preceding paragraph can be chosen so that a fixed Euclidean
profile ball $D$ is mapped inside $Q(q)$, and its smaller concentric
ball contains the inverse image of $B_d(q,2\varepsilon)$ whenever
$0<\varepsilon\leq\varepsilon_*$, for a fixed
$\varepsilon_*<a/8$.  On the latter inverse image the coordinate
distance to $o$ is at most $C\varepsilon$.  The mean-value and
Cauchy inequalities on $D$, applied to $h\circ G_q$, therefore
give, for every $h\in A^2(\Omega)$,
\begin{align}
 \sup_{B_d(q,2\varepsilon)}|h|
       +|h(q)|&\leq Cj_q^{-1/2}\|h\|_{L^2(Q(q))},
       \label{eq:pz-cauchy-value}\\
 \sup_{B_d(q,2\varepsilon)}|h-h(q)|
       &\leq C\varepsilon j_q^{-1/2}
                          \|h\|_{L^2(Q(q))}.
       \label{eq:pz-cauchy-difference}
\end{align}
The constants in these inequalities are independent of
$\varepsilon$.  The factor $j_q^{-1/2}$ is exactly the change of
variables in the $L^2$ norm; the derivative estimate is taken on
fixed profile balls, so it introduces no inverse power of
$\varepsilon$.

Fix $0<\varepsilon\leq\varepsilon_*$.  Choose a maximal
$\varepsilon$-separated set of centers $\{q_\nu\}$ for $d$.
Properness and separability allow a countable locally finite
choice.  The balls $B_d(q_\nu,\varepsilon)$ cover $\Omega$, while
$B_d(q_\nu,\varepsilon/3)$ are pairwise disjoint.  Set
\[
 P_\nu=B_d(q_\nu,2\varepsilon),\qquad
 Q_\nu=B_d(q_\nu,a),\qquad j_\nu=j_{q_\nu}.
\]
Decrease $a$ and $\varepsilon_*$, if necessary, so that
$a+\varepsilon_*<t_*$.  If $x\in Q_\nu$, the disjoint ball
$B_d(q_\nu,\varepsilon/3)$ lies in
$B_d(x,a+\varepsilon/3)$.  The first estimate in
\eqref{eq:pz-local-volumes} thus gives
\begin{equation}\label{eq:pz-large-overlap}
 \sum_\nu\mathbf1_{Q_\nu}(x)\leq C\varepsilon^{-2n}.
\end{equation}
The same argument with radius $3\varepsilon$ shows that the
$P_\nu$ have an overlap bounded by an absolute constant, independent
of $\varepsilon$.

Choose a nonnegative Lipschitz function $b$ on $[0,\infty)$,
equal to one on $[0,1]$ and zero on $[2,\infty)$, and define
\[
 \beta_\nu(z)=b(d(z,q_\nu)/\varepsilon),\qquad
 D_\varepsilon(z)=\sum_\mu\beta_\mu(z),\qquad
 \chi_\nu=\beta_\nu/D_\varepsilon.
\]
Then $1\leq D_\varepsilon\leq C$, $0\leq\chi_\nu\leq1$,
$\sum_\nu\chi_\nu=1$, and $\operatorname{supp}\chi_\nu
\subset\overline P_\nu$.  Smoothness is unnecessary: these
functions enter only integrals, and no derivative of a partition
function is applied to a holomorphic function or to $f$.

We need one precise consequence of the profile hypothesis for
this partition.  For the fixed $\varepsilon$, the functions
\begin{equation}\label{eq:pz-profile-tests}
 \chi_\nu\circ G_{q_\nu}
\end{equation}
are supported in one common compact subset of the model domain
and form a precompact family in $C^0$, hence in $L^1$, as
$q_\nu$ tends to the boundary.  To see this without making a
regularity assumption on distance, note that $d(\,\cdot\,,q_\mu)$
is $1$-Lipschitz for $d$.  The profile metric upper bound makes
each $\beta_\mu\circ G_{q_\nu}$ Lipschitz with constant
$C/\varepsilon$ on a fixed convex profile neighborhood.  At any
point only a uniformly bounded number of these functions is
nonzero, by the small-ball packing estimate.  When comparing
the denominator at two points, only the union of those two
finite sets of nonzero terms contributes.  Thus the denominator
and the quotient in \eqref{eq:pz-profile-tests} have a uniform
Lipschitz bound for the fixed $\varepsilon$.  The support stays
strictly inside the chosen profile neighborhood, so extension
by zero there preserves equicontinuity.  Arzel\`a--Ascoli proves
the assertion.

It follows that, for each fixed $\varepsilon$,
\begin{equation}\label{eq:pz-small-partition-means}
 b_\nu:=\frac{1}{V(P_\nu)}
                  \int_\Omega f\chi_\nu\,dV\longrightarrow0
          \quad(q_\nu\to\partial\Omega).
\end{equation}
Indeed, if this failed, choose a boundary sequence of centers and
pass to a weak-$*$ profile subsequence, which is zero by hypothesis,
and to an $L^1$-convergent subsequence of
\eqref{eq:pz-profile-tests}.  The constant-Jacobian change of
variables and \eqref{eq:pz-local-volumes} give
\[
 b_\nu=\frac{j_\nu}{V(P_\nu)}
       \int (f\circ G_{q_\nu})(\chi_\nu\circ G_{q_\nu})
                                                 \,dV_{\rm mod},
 \qquad \frac{j_\nu}{V(P_\nu)}\leq C\varepsilon^{-2n}.
\]
The integral tends to zero.  The displayed prefactor is bounded
because $\varepsilon$ is fixed, proving
\eqref{eq:pz-small-partition-means}.

We finish with an estimate for the Toeplitz form.  For $h,k\in
A^2(\Omega)$, equations \eqref{eq:pz-cauchy-value}--
\eqref{eq:pz-cauchy-difference} imply
\begin{equation*}
 \sup_{P_\nu}|h\overline k-h(q_\nu)\overline{k(q_\nu)}|
 \leq C\varepsilon j_\nu^{-1}
                \|h\|_{L^2(Q_\nu)}\|k\|_{L^2(Q_\nu)}.
\end{equation*}
Since $V(P_\nu)\leq Cj_\nu\varepsilon^{2n}$, summing the
freezing error over any subset $I$ of the indices gives
\begin{align}
 &\sum_{\nu\in I}\left|\int f\chi_\nu
       \bigl(h\overline k-h(q_\nu)\overline{k(q_\nu)}\bigr)dV
                                                 \right|\nonumber\\
 &\qquad\leq C\varepsilon^{2n+1}\|f\|_\infty
             \sum_\nu\|h\|_{L^2(Q_\nu)}\|k\|_{L^2(Q_\nu)}
       \leq C\varepsilon\|f\|_\infty\|h\|_2\|k\|_2.
       \label{eq:pz-total-freezing}
\end{align}
The last inequality is Cauchy--Schwarz followed by
\eqref{eq:pz-large-overlap}.  In particular its constant is
independent of $\varepsilon$.  The constant part is bounded by
\begin{align}
 &\sum_{\nu\in I}|h(q_\nu)k(q_\nu)|
                       \left|\int f\chi_\nu\,dV\right|\nonumber\\
 &\qquad\leq C\varepsilon^{2n}
       \sup_{\nu\in I}|b_\nu|
       \sum_\nu\|h\|_{L^2(Q_\nu)}\|k\|_{L^2(Q_\nu)}
       \leq C\sup_{\nu\in I}|b_\nu|\,\|h\|_2\|k\|_2.
       \label{eq:pz-total-constant}
\end{align}
Again the constant is independent of $\varepsilon$.  These estimates
also justify absolute convergence of all the sums used in the form.

For $t>0$, there are only finitely many centers with
$\delta(q_\nu)\geq t$: they form an $\varepsilon$-separated
subset of a compact metric set.  Their partition functions have
compact support in $\Omega$, by properness.  Hence
\[
 K_{t,\varepsilon}:=
       T_{\,f\sum_{\delta(q_\nu)\geq t}\chi_\nu}
\]
is compact.  Indeed restriction of the unit ball of $A^2(\Omega)$
to a compact interior set is relatively compact in $L^2$ by the
local mean-value inequality and Montel's theorem; multiplication
by a bounded compactly supported symbol and then $P_\Omega$
preserve compactness.
Apply \eqref{eq:pz-total-freezing} and
\eqref{eq:pz-total-constant} to the remaining indices.  Taking
the supremum over unit vectors $h,k$ yields
\[
 \|T_f-K_{t,\varepsilon}\|
 \leq C\varepsilon\|f\|_\infty
        +C\sup_{\delta(q_\nu)<t}|b_\nu|.
\]
For each fixed $\varepsilon$, the second term tends to zero as
$t\downarrow0$ by \eqref{eq:pz-small-partition-means}.  Consequently
\[
 \|T_f\|_{\mathrm{ess}}\leq C\varepsilon\|f\|_\infty
       \qquad(0<\varepsilon\leq\varepsilon_*),
\]
with $C$ independent of $\varepsilon$.  Letting
$\varepsilon\downarrow0$ proves compactness.
\end{proof}

\begin{proof}[Proof of Theorem~\ref{thm:psc-fixed-radius}]
If $T_f$ is compact, weak convergence of normalized kernels in
Lemma~\ref{lem:psc-profile-berezin} gives
$\widetilde f\to0$.  That lemma identifies this condition with
zero boundary profiles; Lemma~\ref{lem:psc-profile-means}
identifies it with the fixed-radius mean condition.
Conversely, zero profiles imply compactness by
Lemma~\ref{lem:profile-zero-compactness}.  This proves all the
equivalences and the independence of the chosen positive radius.
\end{proof}

\Needspace{12\baselineskip}
\begin{corollary}[A canonical representative when one mean has vanishing oscillation]
\label{cor:psc-vo-representative}
For $f\in L^\infty(\Omega)$, suppose
$a=\mathcal A_R^\Omega f$ has vanishing oscillation at the boundary:
for every $t>0$,
\[
 \sup_{w\in D_t(z)}|a(w)-a(z)|\longrightarrow0
          \quad(z\to\partial\Omega).
\]
Then $a$ is bounded and continuous, $T_f-T_a$ is compact, and,
for every $s>0$,
\[
 \mathcal A_s^\Omega f-a\in C_0(\Omega).
\]
Thus all positive-radius means determine the same representative
modulo $C_0(\Omega)$ under this oscillation hypothesis.
\end{corollary}

\begin{proof}
For completeness, $\mathcal A_t^\Omega$ sends bounded measurable
functions to continuous functions.  The metric is complete and
proper by the geometry established above.  Its distance sphere
of radius $t$ about an interior point is contained in the image,
under the smooth exponential map, of the radius-$t$ sphere in
the tangent space; it therefore has Euclidean measure zero.
If $z_j\to z$, the triangle inequality sandwiches their radius-$t$
balls between the balls about $z$ with radii
$t-d(z_j,z)$ and $t+d(z_j,z)$.  Their indicators consequently
converge in $L^1(dV_{\rm Euc})$.  Both the numerator and the
positive denominator defining the mean are continuous.
Also $|\mathcal A_t^\Omega f|\leq\|f\|_\infty$.

Vanishing oscillation implies
$\mathcal A_t^\Omega a-a\in C_0(\Omega)$ for every $t>0$.
In particular
\[
 \mathcal A_R^\Omega(f-a)=a-\mathcal A_R^\Omega a
            \longrightarrow0.
\]
Theorem~\ref{thm:psc-fixed-radius} gives compactness of
$T_f-T_a$.  The same theorem, now at radius $s$, gives
$\mathcal A_s^\Omega(f-a)\to0$.  Adding
$\mathcal A_s^\Omega a-a\to0$ proves the assertion; continuity
follows from the first paragraph.
\end{proof}

\section{The polydisc}
\label{pd:scope}

The hypothesis of strong pseudoconvexity cannot simply be removed from
the fixed-radius criterion.  We show this by considering the polydisc.
The necessity of the criterion survives, but its sufficiency fails for
real bounded smooth symbols, already on the bidisc.  The averaging sets
in this section are balls for the Bergman geodesic distance, rather than
products of one-dimensional balls.  In contrast, we prove that one
product box with arbitrary fixed positive radii does characterize
compactness.  For geodesic balls, a sequence of radii tending to zero
also gives an exact criterion.

We use Euclidean Lebesgue measure $dV$ on $\mathbb D^n$ and write $d_B$
for the distance in~\eqref{eq:metric-convention}.  Thus
\[
 ds_B^2=\frac2{n+1}\sum_{\ell=1}^n
       \frac{|dz_\ell|^2}{(1-|z_\ell|^2)^2}.
\]
For $R>0$ set
\begin{equation}\label{pd:mean}
 \begin{gathered}
B_R^B(a)=\{z\in\mathbb D^n:d_B(z,a)<R\},
 \\
 A_R f(a)=\frac{1}{V(B_R^B(a))}
             \int_{B_R^B(a)}f\,dV.
 \end{gathered}
\end{equation}
Multiplying the metric by a fixed constant only rescales the radii in
the statements below.  The limit $a\to\partial\mathbb D^n$ always means
approach to the entire topological boundary, and not merely to the
distinguished boundary.

\begin{proposition}\label{pd:necessary}
Let $f\in L^\infty(\mathbb D^n)$.  If $T_f$ is compact on
$A^2(\mathbb D^n)$, then, for every $R>0$,
\[
 A_R f(a)\longrightarrow0
 \qquad(a\to\partial\mathbb D^n).
\]
\end{proposition}

\begin{proof}
For $a\in\mathbb D^n$, let
\[
 \varphi_a(z)_\ell=\frac{a_\ell-z_\ell}
                         {1-\overline{a_\ell}z_\ell},
 \qquad
 U_a h=(h\circ\varphi_a)J_{\mathbb C}\varphi_a.
\]
The operators $U_a$ are unitary involutions on $A^2(\mathbb D^n)$ and
\[
 U_a^*T_fU_a=T_{f\circ\varphi_a}.
\]
For every fixed $h\in A^2(\mathbb D^n)$, one has $U_a h\rightharpoonup0$
as $a\to\partial\mathbb D^n$.  Indeed, for polynomial $h$ this follows
from convergence to zero uniformly on compact subsets, using
\[
 |J_{\mathbb C}\varphi_a(z)|
 =\prod_{\ell=1}^n
       \frac{1-|a_\ell|^2}{|1-\overline{a_\ell}z_\ell|^2}.
\]
At least one factor tends to zero along every subsequence on which
$a$ has a boundary limit.  The general case follows by polynomial
approximation and unitarity.

It follows from compactness that, for fixed $h,k\in A^2(\mathbb D^n)$,
\[
 \langle T_{f\circ\varphi_a}h,k\rangle
 =\langle T_fU_a h,U_a k\rangle\longrightarrow0.
\]
Consequently every weak-star cluster limit $F$ of
$f\circ\varphi_a$ satisfies
$\int Fz^\alpha\overline{z^\beta}\,dV=0$ for all multiindices
$\alpha,\beta$.  Polynomials in $z$ and $\bar z$ are dense in
$L^1(\mathbb D^n)$, so $F=0$.  We have proved
\begin{equation}\label{pd:profiles}
 f\circ\varphi_a\stackrel{*}{\rightharpoonup}0
 \quad\hbox{in }L^\infty(\mathbb D^n)
 \qquad(a\to\partial\mathbb D^n).
\end{equation}

Put $E_R=B_R^B(0)$.  The automorphisms $\varphi_a$ preserve $d_B$, and
the factors $\prod_\ell(1-|a_\ell|^2)^2$ cancel in the quotient defining
the mean.  Thus
\[
 \begin{gathered}
A_R f(a)=\int_{E_R}(f\circ\varphi_a)(w)q_a(w)\,dV(w),
 \\
 q_a(w)=\frac{\prod_\ell|1-\overline{a_\ell}w_\ell|^{-4}}
 {\displaystyle\int_{E_R}
       \prod_\ell|1-\overline{a_\ell}v_\ell|^{-4}\,dV(v)}.
 \end{gathered}
\]
Since $E_R\Subset\mathbb D^n$, this family of densities extends
continuously in $L^1(E_R)$ to $a\in\overline{\mathbb D^n}$.
Along an arbitrary boundary sequence, pass to a subsequence on which
$a$ converges.  Equation~\eqref{pd:profiles} and this $L^1$ continuity
give the asserted limit.  The subsequence argument proves the limit
along the whole boundary.
\end{proof}

\begin{lemma}\label{pd:box-injectivity}
For $\alpha\in\overline{\mathbb D}$ put
$q_\alpha(w)=|1-\overline\alpha w|^{-4}$.  Fix $r>0$ and let
$D_B(b,r)$ denote a one-dimensional Bergman geodesic disc.  If
$g\in L^\infty(\mathbb D)$ and
\[
 \int_{D_B(b,r)}g(w)q_\alpha(w)\,dA(w)=0
 \qquad(b\in\mathbb D),
\]
then $g=0$ almost everywhere.  Consequently, for
$\alpha\in\overline{\mathbb D^n}$ and fixed $r_1,\ldots,r_n>0$,
the identities
\[
 \int_{\prod_{\ell=1}^nD_B(b_\ell,r_\ell)}
 F(w)\prod_{\ell=1}^nq_{\alpha_\ell}(w_\ell)\,dV(w)=0
 \qquad(b\in\mathbb D^n)
\]
imply $F=0$ for every $F\in L^\infty(\mathbb D^n)$.
\end{lemma}

\begin{proof}
For the first assertion, use invariant area
$d\lambda(w)=(1-|w|^2)^{-2}dA(w)$ and put
\[
 H(w)=(1-|w|^2)^2q_\alpha(w)g(w).
\]
The hypothesis says that the invariant mean of $H$ over every disc
of radius $r$ is zero.  If $|\alpha|=1$, then
\[
 H=P_\alpha^2g,
 \qquad P_\alpha(w)=\frac{1-|w|^2}{|1-\overline\alpha w|^2},
\]
so $|H|\le\|g\|_\infty P_\alpha^2$.
If $|\alpha|<1$, the function $q_\alpha$ is bounded.  For any fixed
$\zeta\in\partial\mathbb D$, the inequality
$|1-\overline\zeta w|\le2$ gives
\[
 (1-|w|^2)^2\le16P_\zeta(w)^2.
\]
Thus $|H|\le C P_\zeta^2$ in this case as well.  Apply
Theorem~\ref{sp:one-radius} in complex dimension one, with $m=2$;
the harmless constant difference between the metric normalizations
only changes the radius to another positive radius.  Since
$2>1/2$, that theorem gives $H=0$, and hence $g=0$.

We prove the product assertion by induction.  Fix all centers except
$b_1$ and integrate $F$ in $w_2,\ldots,w_n$ over the corresponding
discs, including their weights.  This gives a bounded function of
$w_1$: all the integration sets are relatively compact and all
weights are bounded on them.  The first assertion shows that this
function is zero for almost every $w_1$.

To make the exceptional set independent of the remaining centers,
first take those centers in a fixed countable dense subset of
$\mathbb D^{n-1}$, and discard the union of the resulting null sets.
For each remaining $w_1$, the weighted integrals depend continuously
on the centers.  Indeed, locally in the center parameters the
integration sets lie in a common compact set, the weights are
bounded, and their indicator functions converge almost everywhere;
the boundaries of these discs have area zero.  Dominated convergence
extends the identities to all remaining centers.  The induction
hypothesis now gives $F(w_1,\cdot)=0$ almost everywhere.  Fubini's
theorem completes the proof.
\end{proof}

For $\mathbf r=(r_1,\ldots,r_n)\in(0,\infty)^n$, define the product
box and its Euclidean mean by
\begin{equation*}
 Q_{\mathbf r}(a)=\prod_{\ell=1}^nD_B(a_\ell,r_\ell),
 \qquad
 A_{\mathbf r}^{\mathrm{box}}f(a)
 =\frac{1}{V(Q_{\mathbf r}(a))}
       \int_{Q_{\mathbf r}(a)}f\,dV.
\end{equation*}

\begin{theorem}\label{pd:box}
Let $f\in L^\infty(\mathbb D^n)$ and fix arbitrary positive radii
$r_1,\ldots,r_n$.  Then
\[
 T_f\text{ is compact on }A^2(\mathbb D^n)
 \quad\Longleftrightarrow\quad
 A_{\mathbf r}^{\mathrm{box}}f(a)\longrightarrow0
 \quad(a\to\partial\mathbb D^n).
\]
\end{theorem}

\begin{proof}
For necessity, the proof of Proposition~\ref{pd:necessary} applies
with $E_R$ replaced by the fixed relatively compact product box
$Q_{\mathbf r}(0)$.  Its pulled-back normalized densities are again
continuous in $L^1$ for $a\in\overline{\mathbb D^n}$.

For sufficiency, take an arbitrary sequence $a^{(\nu)}$ approaching
the boundary.  After passing to a subsequence, suppose that
\[
 a^{(\nu)}\longrightarrow\alpha\in\partial\mathbb D^n,
 \qquad
 f\circ\varphi_{a^{(\nu)}}\stackrel{*}{\rightharpoonup}F
 \quad\hbox{in }L^\infty(\mathbb D^n).
\]
For every fixed $b\in\mathbb D^n$, the points
$\varphi_{a^{(\nu)}}(b)$ also approach the boundary.  This follows
coordinatewise from
\[
 1-|\varphi_{a_\ell}(b_\ell)|^2
 =\frac{(1-|a_\ell|^2)(1-|b_\ell|^2)}
        {|1-\overline{a_\ell}b_\ell|^2}
\]
in any coordinate for which $|\alpha_\ell|=1$.
Pull the hypothesized mean back by $\varphi_{a^{(\nu)}}$.  Its domain
becomes $Q_{\mathbf r}(b)$, and its constant Jacobian factors cancel
between numerator and denominator.  The remaining densities
converge uniformly on this compact integration set to
$\prod_\ell q_{\alpha_\ell}$.  Their integrals converge to a strictly
positive finite denominator.  The weak-star convergence therefore
gives
\[
 \int_{Q_{\mathbf r}(b)}
 F(w)\prod_{\ell=1}^nq_{\alpha_\ell}(w_\ell)\,dV(w)=0
 \qquad(b\in\mathbb D^n).
\]
Lemma~\ref{pd:box-injectivity} yields $F=0$.

Thus every boundary automorphic profile is zero.  In particular,
using the constant normalized kernel at the origin, every such
subsequence satisfies
\[
 \widetilde f(a^{(\nu)})
 =\pi^{-n}\int_{\mathbb D^n}
        (f\circ\varphi_{a^{(\nu)}})(w)\,dV(w)
 \longrightarrow0.
\]
The subsequence argument shows $\widetilde f(a)\to0$ along the whole
boundary.  Since $T_f$ belongs to the Toeplitz algebra on the
polydisc, \cite[Theorem~1.1]{MW2014} implies that $T_f$ is compact.
\end{proof}

\begin{corollary}\label{pd:shrinking-radii}
Let $f\in L^\infty(\mathbb D^n)$ and let $(R_k)$ be any positive
sequence with $R_k\downarrow0$.  Then
\[
 T_f\text{ is compact}
 \quad\Longleftrightarrow\quad
 \text{for every }k,\quad
 A_{R_k}f(a)\longrightarrow0
 \quad(a\to\partial\mathbb D^n).
\]
Here $A_{R_k}$ uses geodesic balls as in~\eqref{pd:mean}.  The
boundary limit is taken separately for each fixed $k$.
\end{corollary}

\begin{proof}
Necessity follows from Proposition~\ref{pd:necessary}.  For
sufficiency, take $\alpha$ and a profile $F$ along a boundary
subsequence as in the preceding proof.  The same change of variables
and limit argument gives
\[
 \int_{B_{R_k}^B(b)}F(w)
       \prod_{\ell=1}^nq_{\alpha_\ell}(w_\ell)\,dV(w)=0
 \qquad(b\in\mathbb D^n,\ k\ge1).
\]
The function
$F\prod_\ell q_{\alpha_\ell}$ is locally integrable.  Smooth
Riemannian balls are a differentiation basis for Euclidean
Lebesgue measure in the interior: on each compact coordinate
neighborhood, sufficiently small metric balls are contained in
Euclidean balls of comparable radii and have comparable volumes.
Lebesgue differentiation along $R_k\downarrow0$ consequently gives
$F\prod_\ell q_{\alpha_\ell}=0$ almost everywhere.  All these weights
are positive in the interior, so $F=0$.  The Berezin-transform
argument in Theorem~\ref{pd:box}, followed by
\cite[Theorem~1.1]{MW2014}, proves compactness.
\end{proof}

\begin{corollary}[Nonnegative symbols]\label{pd:positive}
Let $f\in L^\infty(\mathbb D^n)$ be nonnegative and fix $R>0$.
Then $T_f$ is compact if and only if $A_Rf(a)\to0$ as
$a\to\partial\mathbb D^n$.
\end{corollary}
\begin{proof}
Necessity is Proposition~\ref{pd:necessary}. Choose a fixed positive
box radius vector $\mathbf r$ so that
$Q_{\mathbf r}(a)\subset B_R^B(a)$ for every $a$.
The quotient $V(B_R^B(a))/V(Q_{\mathbf r}(a))$ is bounded uniformly
in $a$: pull both sets back by $\varphi_a$, cancel the common
Jacobian factor, and use the positive, locally uniform limiting
weights for $a\in\overline{\mathbb D^n}$ on these fixed compact sets.
Positivity gives $A_{\mathbf r}^{\mathrm{box}}f(a)\le C A_Rf(a)$.
Theorem~\ref{pd:box} now proves sufficiency.
\end{proof}

\Needspace{18\baselineskip}
\begin{theorem}\label{pd:counterexample}
For every $n\ge2$ there exists $R_n>0$ such that, for each fixed
$0<R<R_n$, there are a number $t=t_R>0$ and a real-valued function
\[
 f\in C^\infty(\mathbb D^n)\cap L^\infty(\mathbb D^n),
 \qquad |f|\le1,
\]
and a sequence $a_j\to\partial\mathbb D^n$ for which
\[
 \begin{gathered}
A_R f(a)\longrightarrow0
 \quad(a\to\partial\mathbb D^n),
 \\
 \widetilde f(a_j)\longrightarrow
 c_t:=(1+t^2)\left(\frac{\pi t}{\sinh(\pi t)}\right)^2>0.
 \end{gathered}
\]
In particular,
\[
 \|T_f\|_{\mathrm{ess}}\ge c_t>0.
\]
Here $A_R$ is the Euclidean-volume average over the true Bergman
geodesic ball in~\eqref{pd:mean}, the boundary limit is taken along
the entire topological boundary, and the asserted smoothness is
inside $\mathbb D^n$.  In particular, these averages do not
characterize compactness on the polydisc for general real or
complex bounded symbols.
\end{theorem}

\begin{proof}
We first construct a null function for the corresponding average on
a product of upper half-planes.  Write $\mathbb H=\{x+iy:y>0\}$ and
let $d_0$ be the product distance associated with
\[
 ds_0^2=\sum_{\ell=1}^n
                  \frac{dx_\ell^2+dy_\ell^2}{y_\ell^2}.
\]
The Bergman metric on $\mathbb H^n$, in our convention, is
$ds_B^2=\tfrac1{2(n+1)}ds_0^2$.  Hence the Cayley transform identifies a
Bergman ball of radius $R$ with a $d_0$-ball of radius
$r=\sqrt{2(n+1)}R$.

Let $\boldsymbol i=(i,\ldots,i)$ and, for real $t$, define
\begin{equation*}
 F_t(z)=\left(\frac{y_1}{y_2}\right)^{it},
 \qquad
 I_r(t)=\int_{B_0(\boldsymbol i,r)}F_t(z)\,dV(z).
\end{equation*}
Here $B_0$ denotes a ball for $d_0$.  Interchanging the first two
coordinates conjugates the integrand and preserves the integration
set and measure.  Therefore $I_r(t)$ is real.  It is continuous in
$t$, and $I_r(0)>0$.

For fixed $u\in\mathbb R$, use coordinates
\[
 z_\ell=rx_\ell+i(1+rv_\ell).
\]
The rescaled balls tend to the unit ball in $\mathbb R^{2n}$, while
\[
 F_{u/r}(z)
 =\exp\left\{\frac{iu}{r}
       [\log(1+rv_1)-\log(1+rv_2)]\right\}
 \longrightarrow e^{iu(v_1-v_2)}.
\]
These convergences hold uniformly on fixed compact sets, apart from
the irrelevant indicator discontinuity on the limiting sphere.
For small $r$, the rescaled balls lie in a fixed bounded set:
indeed $d_{\mathbb H}(i,z_\ell)<r$ implies
$|z_\ell-i|\le e^r-1$.  Dominated convergence consequently gives
\begin{equation}\label{pd:Bessel}
 r^{-2n}I_r(u/r)\longrightarrow
 \int_{|(x,v)|<1}e^{iu(v_1-v_2)}\,dx\,dv
 =(2\pi)^n\frac{J_n(\sqrt2u)}{(\sqrt2u)^n},
\end{equation}
with the value at $u=0$ understood by continuity.  The standard
Fourier transform formula for a Euclidean ball gives the last
identity.

Choose $u_*>0$ for which $J_n(\sqrt2u_*)<0$.  Equation~\eqref{pd:Bessel}
shows that there exists $r_n>0$ such that
$I_r(u_*/r)<0$ whenever $0<r<r_n$.  Since $I_r(0)>0$, the intermediate
value theorem supplies $t=t_r>0$ with $I_r(t)=0$.
Set $t_R=t_r$ for the corresponding Bergman radius
$R=r/\sqrt{2(n+1)}$.

For $b=(\alpha_\ell+i\beta_\ell)_{\ell=1}^n\in\mathbb H^n$, the map
\[
 L_b(z)_\ell=\alpha_\ell+\beta_\ell z_\ell
\]
is a $d_0$-isometry taking $\boldsymbol i$ to $b$.  Its real Jacobian
is $\prod_\ell\beta_\ell^2$, and
$F_t\circ L_b=(\beta_1/\beta_2)^{it}F_t$.  Thus
\begin{equation}\label{pd:modelzero}
 \int_{B_0(b,r)}F_t\,dV=0
 \qquad\text{for every }b\in\mathbb H^n.
\end{equation}
The real function $G=\operatorname{Re}F_t$ satisfies the same identity.

We also record that the Berezin transform of $G$ does not vanish at
$\boldsymbol i$.  On $\mathbb H$ the normalized Bergman kernel at $i$
has density $4/(\pi|z+i|^4)$.  Integrating first in $x$ gives the
height density $2/(1+y)^3$, so
\begin{equation*}
 c(t)=2\int_0^\infty\frac{y^{it}}{(1+y)^3}\,dy
     =\Gamma(1+it)\Gamma(2-it)
     =(1-it)\frac{\pi t}{\sinh(\pi t)}.
\end{equation*}
The kernel on the product is the product kernel.  It follows that
\begin{equation}\label{pd:positiveBerezin}
 \widetilde G(\boldsymbol i)=c(t)c(-t)
     =c_t=(1+t^2)\left(\frac{\pi t}{\sinh(\pi t)}\right)^2>0.
\end{equation}

It remains to place this model in $\mathbb D^n$.  This step must
account for the nonconstant Euclidean Jacobian of the Cayley map.
Let
\[
 \kappa(w)=\frac{i-w}{i+w},\qquad
 \mathcal C=\kappa^{\times n},\qquad
 J(w)=|J_{\mathbb C}\mathcal C(w)|^2
     =\prod_{\ell=1}^n\frac4{|i+w_\ell|^4}.
\]
The function $J$ extends smoothly and positively to a neighborhood
of $0\in\mathbb C^n$, with $J(0)=4^n$.

Choose $\eta\in C^\infty([0,\infty))$ satisfying
$0\le\eta\le1$, $\eta=1$ on $[0,1]$, and $\eta=0$ on $[2,\infty)$.
For integers $j$ sufficiently large put
\begin{equation*}
 \begin{gathered}
L_j=j,\qquad \delta_j=e^{-j^3},\qquad
 p_j=(i\delta_j,\ldots,i\delta_j),\\
 \chi_j(w)=\eta\!\left(\frac{d_0(w,p_j)}{L_j}\right).
 \end{gathered}
\end{equation*}
These cutoffs are smooth, since the product hyperbolic space is a
Hadamard manifold, its distance from a point is smooth away from
that point, and $\eta$ is constant near zero.  They are Lipschitz
with constants at most $C/L_j$.  Their supports are contained in
$B_0(p_j,2L_j)$.  Since
\[
 d_0(p_j,p_k)=\sqrt n\,|j^3-k^3|,
\]
we may discard finitely many indices so that these supporting balls
are pairwise separated by more than $4r+1$.

If $w\in B_0(p_j,2L_j+2r)$, then each coordinate satisfies
\[
 |w_\ell|\le\delta_j e^{2L_j+2r}=:\varepsilon_j,
 \qquad \varepsilon_j\longrightarrow0.
\]
Indeed, the hyperbolic ball of radius $s$ about $i\delta_j$ lies in
$\{|w|\le\delta_je^s\}$.  In particular, the supports escape every
compact subset of $\mathbb H^n$ and their images approach the single
boundary point $(1,\ldots,1)$.

Define
\[
 g(w)=\sum_j\chi_j(w)G(w),\qquad
 f=g\circ\mathcal C^{-1}.
\]
The sum is locally finite and its summands have disjoint supports.
Thus $f$ is real, belongs to $C^\infty(\mathbb D^n)$, and satisfies
$|f|\le1$.

We next check the mean along the entire boundary.  A ball
$B_0(b,r)$ meets at most one supporting ball.  If it meets the
$j$th one, it is contained in $B_0(p_j,2L_j+2r)$, and
\[
 A_R f(\mathcal C(b))
 =\frac{\displaystyle\int_{B_0(b,r)}\chi_j(w)G(w)J(w)\,dV(w)}
        {\displaystyle\int_{B_0(b,r)}J(w)\,dV(w)}.
\]
By~\eqref{pd:modelzero}, we may subtract
$\chi_j(b)J(b)G(w)$ inside the numerator.  The Lipschitz bound for
$\chi_j$ and smoothness of $J$ at zero then imply
\begin{equation}\label{pd:meanestimate}
 |A_R f(\mathcal C(b))|
 \le C\left(\frac r{L_j}+\varepsilon_j\right).
\end{equation}
Here the constant is independent of $j$ and $b$; $r$ is fixed, and
$J$ is bounded above and below by positive constants on the
relevant neighborhoods.  A ball meeting none of the supports has
mean zero.  Along any sequence $\mathcal C(b)\to\partial\mathbb D^n$,
the indices of the supports met must tend to infinity, since a
fixed finite collection has a compact $r$-neighborhood.  Hence
\eqref{pd:meanestimate} proves $A_R f\to0$ along the whole boundary.

Finally set $a_j=\mathcal C(p_j)$.  Common dilation of all coordinates
is a biholomorphism of $\mathbb H^n$, and
\[
 \chi_j(\delta_jw)
 =\eta\!\left(\frac{d_0(w,\boldsymbol i)}{L_j}\right),
 \qquad G(\delta_jw)=G(w).
\]
Since $L_j=j$, the first identity gives
$\chi_j(\delta_jw)=1$ throughout $B_0(\boldsymbol i,j)$.
Disjointness of the supporting balls consequently implies
\[
 g(\delta_jw)=G(w)
 \qquad\bigl(w\in B_0(\boldsymbol i,j)\bigr).
\]
Introduce the probability measure
\[
 d\nu(w)=
 \prod_{\ell=1}^n\frac4{\pi|w_\ell+i|^4}\,dV(w),
\]
which is the normalized Bergman kernel measure at
$\boldsymbol i$.  Berezin transforms are invariant under
biholomorphic changes of coordinates, including the accompanying
Jacobian in the unitary map on Bergman spaces.  Hence
\[
 \widetilde f(a_j)=\widetilde g(p_j)
      =\int_{\mathbb H^n}g(\delta_jw)\,d\nu(w),
 \qquad
 c_t=\int_{\mathbb H^n}G(w)\,d\nu(w).
\]
As $|g|,|G|\le1$, it follows that
\[
 |\widetilde f(a_j)-c_t|
 \le 2\nu\!\left(\mathbb H^n\setminus
                         B_0(\boldsymbol i,j)\right).
\]
We estimate this tail explicitly.  In each coordinate, the Cayley
map $\kappa$ sends the corresponding one-dimensional factor of
$\nu$ to normalized area measure on $\mathbb D$ and satisfies
\[
 d_{\mathbb H}(i,w)=2\operatorname{arctanh}|\kappa(w)|.
\]
Consequently, for $s\ge0$, that factor assigns mass
\[
 1-\tanh^2(s/2)=\frac1{\cosh^2(s/2)}
\]
to $\{w:d_{\mathbb H}(i,w)\ge s\}$.
The product-distance formula shows that a point outside
$B_0(\boldsymbol i,j)$ has at least one coordinate at hyperbolic
distance at least $j/\sqrt n$ from $i$.  The union bound therefore
gives
\[
 |\widetilde f(a_j)-c_t|
 \le \frac{2n}{\cosh^2\!\bigl(j/(2\sqrt n)\bigr)}
 \le 8n e^{-j/\sqrt n}.
\]
In particular, $\widetilde f(a_j)\to c_t>0$.

The normalized kernels $k_{a_j}$ converge weakly to zero on
$A^2(\mathbb D^n)$, as follows from their explicit product formula,
unit norm, and convergence to zero on compact subsets.  If $K$ is
any compact operator on this space, then
$\langle Kk_{a_j},k_{a_j}\rangle\to0$.  Thus
\[
 \|T_f-K\|
 \ge\limsup_{j\to\infty}
       \left|\langle (T_f-K)k_{a_j},k_{a_j}\rangle\right|
 =c_t.
\]
Taking the infimum over compact $K$ proves
$\|T_f\|_{\mathrm{ess}}\ge c_t$.
Taking $R_n=r_n/\sqrt{2(n+1)}$ completes the proof.
\end{proof}

Cancellation between two coordinate directions produces the null
function in~\eqref{pd:modelzero}. A geodesic ball couples the
coordinate distances through their squared sum. Its small-scale
multiplier is the sign-changing Bessel function
in~\eqref{pd:Bessel}, while the Berezin multiplier
in~\eqref{pd:positiveBerezin} is positive.

\section{The weighted unit ball}
\label{sec:weighted-ball}

We prove the weighted version of the fixed-radius criterion on the
unit ball. Let
\[
 dv=V(\mathbb B_n)^{-1}dV,\qquad
 d\lambda(z)=(1-|z|^2)^{-n-1}dv(z),
\]
and put
\[
 \begin{gathered}
dv_\alpha(z)=c_\alpha(1-|z|^2)^\alpha dv(z),\quad
 \alpha>-1,\quad v_\alpha(\mathbb B_n)=1,
 \\ A^2_\alpha= L^2(dv_\alpha)\cap\mathcal O(\mathbb B_n).
 \end{gathered}
\]
The averages in this section use $dv_\alpha$, matching the Hilbert
space measure. Here $D(a,r)=D_{\mathbb B_n}(a,r)$ for the fixed
metric normalization in~\eqref{eq:metric-convention}. Set
\begin{equation}\label{eq:defaverage}
 \mathcal A_{r,\alpha}\phi(a)
 =\frac{1}{v_\alpha(D(a,r))}
       \int_{D(a,r)}\phi\,dv_\alpha.
\end{equation}
\begin{theorem}\label{thm:main}
Let $n\ge1$, $\alpha>-1$, and $\phi\in L^\infty(\mathbb B_n)$.
For every fixed $r>0$, $T_\phi$ is compact on $A^2_\alpha$ if and
only if
\begin{equation}\label{eq:averagevanish}
 \mathcal A_{r,\alpha}\phi(a)\longrightarrow0\qquad(|a|\to1).
\end{equation}
\end{theorem}
\begin{corollary}\label{cor:radius}
The condition in~\eqref{eq:averagevanish} for one positive radius
is equivalent to the same condition for every positive radius.
\end{corollary}
For $n=1$ and $\alpha=0$ the conclusion is contained in
Cao--He--Zhang~\cite{CHZ}. The proof here also covers the other
dimensions and all standard weights.

\subsection{Automorphisms and the Berezin criterion}

We collect the standard automorphism identities needed below.  For $a\in\B_n$, let $\varphi_a$ be the involutive automorphism of $\B_n$ interchanging $0$ and $a$.  With $P_a$ the orthogonal projection onto $\C a$, $Q_a=I-P_a$, and $s_a=(1-|a|^2)^{1/2}$, one may take
\[
\varphi_a(z)=\frac{a-P_a z-s_aQ_a z}{1-\ip{z}{a}}.
\]
Then
\begin{equation}\label{eq:automorphism-identity}
1-|\varphi_a(z)|^2
=\frac{(1-|a|^2)(1-|z|^2)}{|1-\ip{z}{a}|^2},
\end{equation}
and the real Jacobian is
\begin{equation*}
J_{\mathbb R}\varphi_a(z)
=\left(\frac{1-|a|^2}{|1-\ip{z}{a}|^2}\right)^{n+1}.
\end{equation*}
Consequently, with
\[
m=n+1+\alpha,
\]
we have
\begin{equation}\label{eq:weighted-jacobian}
dv_\alpha(\varphi_a(z))
=\left(\frac{1-|a|^2}{|1-\ip{z}{a}|^2}\right)^m dv_\alpha(z).
\end{equation}
See, for example, \cite[Chapter 2]{ZhuBall}.

The normalized reproducing kernel is, up to the harmless normalization already built into $dv_\alpha$,
\begin{equation*}
k_a^\alpha(z)
=\frac{(1-|a|^2)^{m/2}}{(1-\ip{z}{a})^m}.
\end{equation*}
Choose a holomorphic branch of the automorphism multiplier and define
\begin{equation}\label{eq:Ua}
(U_af)(z)
=\frac{(1-|a|^2)^{m/2}}{(1-\ip{z}{a})^m}
 f(\varphi_a(z)).
\end{equation}
Then $U_a$ is unitary on $A^2_\alpha(\B_n)$, $U_a^2$ is multiplication by a unimodular constant, and
\begin{equation}\label{eq:conjugation}
U_aT_\phi U_a^*=T_{\phi\circ\varphi_a}.
\end{equation}
In particular,
\begin{equation}\label{eq:berezin-translate}
\wt\phi_\alpha(a)
=\int_{\B_n}(\phi\circ\varphi_a)(z)\,dv_\alpha(z).
\end{equation}

Let $\mathfrak T_\alpha$ be the norm-closed algebra generated by
$\{T_\phi:\phi\in L^\infty(\mathbb B_n)\}$. We use the following
standard compactness criterion from \cite[Theorem~1.1]{MSW}.

\begin{theorem}[Mitkovski--Su\'arez--Wick]\label{thm:MSW}
Let $\alpha>-1$ and let
$S\in\mathcal B(A^2_\alpha(\B_n))$.  Then $S$ is compact if and only if
$S\in\mathfrak T_\alpha$ and
\[
\widetilde S(a)\longrightarrow0\qquad(|a|\to1).
\]
In particular, for $S\in\mathfrak T_\alpha$, compactness is equivalent
to boundary vanishing of the Berezin transform.
\end{theorem}

\begin{corollary}\label{thm:BerezinCompact}
For $\alpha>-1$ and $\phi\in L^\infty(\B_n)$,
\[
T_\phi\text{ is compact on }A^2_\alpha(\B_n)
\quad\Longleftrightarrow\quad
\widetilde\phi_\alpha(a)\to0\quad(|a|\to1).
\]
\end{corollary}

For later use we also record weak convergence of the automorphism unitaries.

\begin{lemma}\label{lem:Uaweak}
For every $f\in A^2_\alpha(\B_n)$,
\[
U_af\rightharpoonup0\qquad (|a|\to1).
\]
\end{lemma}

\begin{proof}
For a polynomial $p$, the multiplier in \eqref{eq:Ua} tends locally uniformly to zero as $|a|\to1$, while $p\circ\varphi_a$ remains locally bounded.  Thus $U_ap\to0$ locally uniformly.  Testing against reproducing kernels gives weak convergence for polynomials.  Since the $U_a$ are unitary and polynomials are dense, the conclusion follows for arbitrary $f$.
\end{proof}

\subsection{Boundary profiles on the ball}
Throughout this section
\[
m=n+1+\alpha>n.
\]

\begin{proposition}[Boundary annihilation of automorphic translates]\label{prop:weakstar-annihilation}
Let $f\in L^\infty(\B_n)$ and fix $r>0$.  If
\begin{equation}\label{eq:average-vanish-f}
\mathcal A_{r,\alpha}f(a)\longrightarrow0
\qquad(|a|\to1),
\end{equation}
then
\begin{equation}\label{eq:all-translates-weakstar}
f\circ\varphi_a\overset{w^*}{\longrightarrow}0
\qquad\text{in }L^\infty(\B_n,dv_\alpha).
\end{equation}
\end{proposition}

\begin{proof}
Suppose \eqref{eq:all-translates-weakstar} fails.  Since
$L^1(\B_n,dv_\alpha)$ is separable and the translated functions are
uniformly bounded in $L^\infty$, there are $a_j\in\B_n$ with
$|a_j|\to1$ and, after passing to a subsequence,
\[
a_j\to\zeta\in\partial\B_n,
\qquad
f\circ\varphi_{a_j}\overset{w^*}{\longrightarrow}\psi
\]
for some nonzero $\psi\in L^\infty(\B_n)$.  Fix $b\in\B_n$ and put
$c_j=\varphi_{a_j}(b)$.  Then $c_j\to\zeta$, and
$\varphi_{a_j}(D(b,r))=D(c_j,r)$.  Using
\eqref{eq:weighted-jacobian} and cancelling the common factor
$(1-|a_j|^2)^m$, where $m=n+1+\alpha$, we obtain
\[
\mathcal A_{r,\alpha}f(c_j)
=\frac{\displaystyle
\int_{D(b,r)}(f\circ\varphi_{a_j})(z)
|1-\langle z,a_j\rangle|^{-2m}\,dv_\alpha(z)}
{\displaystyle
\int_{D(b,r)}|1-\langle z,a_j\rangle|^{-2m}\,dv_\alpha(z)}.
\]
Since $D(b,r)\Subset\B_n$, the weights converge uniformly there to
$|1-\langle z,\zeta\rangle|^{-2m}$.  Passing to the limit and using
\eqref{eq:average-vanish-f} gives
\[
\int_{D(b,r)}
\frac{\psi(z)(1-|z|^2)^\alpha}
{|1-\langle z,\zeta\rangle|^{2m}}\,dv(z)=0,
\qquad b\in\B_n.
\]
Equivalently, with
\[
F(z)=\psi(z)P_\zeta(z)^m,
\]
we have
\[
\int_{D(b,r)}F\,d\lambda=0
\qquad(b\in\B_n),
\qquad
|F(z)|\le \|f\|_\infty P_\zeta(z)^m.
\]
Because $m=n+1+\alpha>n$, \cref{sp:one-radius} gives $F=0$, and
hence $\psi=0$, a contradiction.  Thus
\eqref{eq:all-translates-weakstar} holds.
\end{proof}

\begin{proposition}\label{prop:sufficiency}
Let $\phi\in L^\infty(\B_n)$.  If \eqref{eq:averagevanish} holds for
one fixed $r>0$, then $T_\phi$ is compact on $A^2_\alpha(\B_n)$.
\end{proposition}

\begin{proof}
By Proposition~\ref{prop:weakstar-annihilation},
$\phi\circ\varphi_a\overset{w^*}{\longrightarrow}0$.  Taking the
constant test function $1$ in \eqref{eq:berezin-translate} gives
$\widetilde\phi_\alpha(a)\to0$.  The conclusion follows from
\cref{thm:BerezinCompact}.
\end{proof}

The converse direction is simpler and does not use complex-hyperbolic rigidity.

\begin{proposition}\label{prop:necessity}
If $T_\phi$ is compact on $A^2_\alpha(\B_n)$, then
\[
\mathcal A_{r,\alpha}\phi(a)\to0
\qquad(|a|\to1)
\]
for every fixed $r>0$.
\end{proposition}

\begin{proof}
By Lemma~\ref{lem:Uaweak}, $U_af\rightharpoonup0$ for every $f\in A^2_\alpha$.  Compactness gives
\[
\|T_\phi U_af\|\to0.
\]
Using \eqref{eq:conjugation}, for $f,g\in A^2_\alpha$,
\begin{align*}
\int_{\B_n}(\phi\circ\varphi_a)f\overline g\,dv_\alpha
&=\ip{T_{\phi\circ\varphi_a}f}{g}\\
&=\ip{T_\phi U_af}{U_ag}\to0.
\end{align*}
The span of $z^\mu\overline z^{\nu}$ is dense in $L^1(\B_n,dv_\alpha)$: by Stone--Weierstrass it is uniformly dense in $C(\overline{\B_n})$, and continuous functions are dense in $L^1(dv_\alpha)$.  Since the translated symbols are uniformly bounded,
\begin{equation}\label{eq:weakstar-zero}
\phi\circ\varphi_a\overset{w^*}{\longrightarrow}0
\qquad\text{in }L^\infty(\B_n).
\end{equation}

For fixed $r>0$, change variables in \eqref{eq:defaverage} using $D(a,r)=\varphi_a(D(0,r))$ and \eqref{eq:weighted-jacobian}:
\begin{equation}\label{eq:average-at-a}
\mathcal A_{r,\alpha}\phi(a)
=\frac{\displaystyle
\int_{D(0,r)}(\phi\circ\varphi_a)(z)|1-\ip{z}{a}|^{-2m}\,dv_\alpha(z)}
{\displaystyle
\int_{D(0,r)}|1-\ip{z}{a}|^{-2m}\,dv_\alpha(z)}.
\end{equation}
If $a_j\to\zeta\in\partial\B_n$, the weight in \eqref{eq:average-at-a} converges uniformly on the compact set $D(0,r)$ to $|1-\ip{z}{\zeta}|^{-2m}$.  Equation \eqref{eq:weakstar-zero} makes the numerator tend to zero, while the denominator tends to a positive finite number.  Hence $\mathcal A_{r,\alpha}\phi(a_j)\to0$.  Since every sequence approaching the boundary has a subsequence converging to some boundary point, the full limit follows.
\end{proof}

\begin{proof}[Proof of \cref{thm:main}]
The implication from the fixed-radius condition to compactness is Proposition~\ref{prop:sufficiency}; the converse is Proposition~\ref{prop:necessity}.
\end{proof}

\begin{proof}[Proof of Corollary~\ref{cor:radius}]
If the condition holds for one radius, Proposition~\ref{prop:sufficiency} gives compactness.  Then Proposition~\ref{prop:necessity} gives the condition for every radius.
\end{proof}

\subsection{Vanishing mass and unimodular symbols}\label{sec:carleson}

For positive symbols, the main theorem is closely related to the
classical theory of vanishing Carleson measures.  Recall that a positive
Borel measure $\mu$ on $\B_n$ is a vanishing Carleson measure for
$A^2_\alpha(\B_n)$ if the embedding
\[
I_\mu:A^2_\alpha(\B_n)\longrightarrow L^2(\mu)
\]
is compact.  Equivalently (see, for example, \cite[Chapter~3]{ZhuBall}),
for one, and hence for every, fixed $r>0$,
\begin{equation}\label{eq:VC-density}
\frac{\mu(D(a,r))}{v_\alpha(D(a,r))}\longrightarrow0
\qquad(|a|\to1).
\end{equation}
If $\phi\ge0$ and $d\mu_\phi=\phi\,dv_\alpha$, then
\[
\frac{\mu_\phi(D(a,r))}{v_\alpha(D(a,r))}
=\mathcal A_{r,\alpha}\phi(a).
\]
Thus, for nonnegative bounded symbols, \cref{thm:main} recovers exactly
the standard vanishing-Carleson characterization of compact positive
Toeplitz operators.

For complex-valued symbols the situation is different.  Condition
\eqref{eq:VC-density} applied to the total variation would require
\begin{equation}\label{eq:VC-absolute}
\mathcal A_{r,\alpha}(|\phi|)(a)\longrightarrow0,
\end{equation}
whereas our theorem only requires
\begin{equation}\label{eq:signed-average}
\mathcal A_{r,\alpha}\phi(a)\longrightarrow0.
\end{equation}
The implication \eqref{eq:VC-absolute}$\Rightarrow$\eqref{eq:signed-average}
is immediate, but the converse fails even for unimodular symbols.

\begin{proposition}[A unimodular symbol in every Schatten class]\label{prop:oscillatory-symbol}
Let $\gamma>0$ and define
\begin{equation}\label{eq:osc-symbol}
\phi_\gamma(z)
=\exp\!\bigl(i(1-|z|^2)^{-\gamma}\bigr),
\qquad z\in\B_n.
\end{equation}
Then $|\phi_\gamma|\equiv1$ and, for every $\alpha>-1$ and every $p>0$,
\begin{equation}\label{eq:all-schatten}
T_{\phi_\gamma}\in\mathcal S_p\bigl(A^2_\alpha(\B_n)\bigr).
\end{equation}
More precisely, if $\lambda_\ell$ is the eigenvalue of
$T_{\phi_\gamma}$ on the space of homogeneous holomorphic polynomials of
degree $\ell$, then for every integer $N\ge1$ there is a constant
$C_N>0$ such that
\begin{equation}\label{eq:eigen-superpoly}
|\lambda_\ell|\le C_N(1+\ell)^{-N\gamma},
\qquad \ell\ge0.
\end{equation}
Consequently the singular values of $T_{\phi_\gamma}$ decay faster than
every power.  On the other hand,
$|\phi_\gamma|\,dv_\alpha=dv_\alpha$ is not a vanishing Carleson measure.
Moreover, for every fixed invariant radius $r>0$,
\begin{equation}\label{eq:osc-average-rate}
\bigl|\mathcal A_{r,\alpha}\phi_\gamma(a)\bigr|
\le C_{r,\gamma}(1-|a|^2)^\gamma
\end{equation}
for $|a|$ sufficiently close to $1$.
\end{proposition}

\begin{proof}
We divide the argument into four steps.  The first one exhibits directly
the fixed-radius cancellation; the remaining steps give the stronger
Schatten conclusion.

\smallskip
\noindent\emph{Step 1: fixed-radius cancellation.}
Fix $r>0$ and put $E=D(0,r)$.  Since $E$ is compactly contained in
$\B_n$, all functions below are smooth on a neighborhood of $\overline E$.
Set
\[
q=n+1+\alpha,\qquad \delta_a=1-|a|^2,
\qquad \Lambda_a=\delta_a^{-\gamma}.
\]
Using $D(a,r)=\varphi_a(E)$, \eqref{eq:weighted-jacobian}, and
\eqref{eq:automorphism-identity}, we obtain
\begin{equation}\label{eq:osc-average-ratio}
\mathcal A_{r,\alpha}\phi_\gamma(a)
=\frac{\displaystyle
\int_E e^{i\Lambda_a Q_a(z)}B_a(z)\,dv_\alpha(z)}
{\displaystyle\int_E B_a(z)\,dv_\alpha(z)},
\end{equation}
where
\begin{equation*}
Q_a(z)=
\left(\frac{|1-\langle z,a\rangle|^2}{1-|z|^2}\right)^\gamma,
\qquad
B_a(z)=|1-\langle z,a\rangle|^{-2q}.
\end{equation*}
The denominator in \eqref{eq:osc-average-ratio} is bounded above and
below by positive constants independent of $a$.

We claim that the numerator is $O(\Lambda_a^{-1})$ as $|a|\to1$.
Indeed, if $a_j\to\zeta\in\partial\B_n$, then on $\overline E$
\[
Q_{a_j}\longrightarrow Q_\zeta
=\left(\frac{|1-\langle z,\zeta\rangle|^2}{1-|z|^2}\right)^\gamma
=P_\zeta(z)^{-\gamma}
\]
in $C^2$, while $B_{a_j}\to B_\zeta$ in $C^1$.  Since
$-\log P_\zeta$ is a Busemann function, its differential is nowhere
zero; hence $\nabla Q_\zeta$ never vanishes.  Compactness of
$\overline E\times\partial\B_n$ and unitary invariance therefore give
constants $c_r>0$ and $a_0<1$ such that
\begin{equation*}
|\nabla Q_a(z)|\ge c_r,
\qquad z\in\overline E,\quad |a|\ge a_0.
\end{equation*}
The relevant $C^2$-norms of $Q_a$ and $C^1$-norms of $B_a$ are also
uniformly bounded.  Absorbing the smooth density of $dv_\alpha$ on $E$
into the amplitude and setting
\[
V_a=\frac{\nabla Q_a}{|\nabla Q_a|^2},
\]
we have $V_a(Q_a)=1$.  One integration by parts gives
\[
\left|\int_E e^{i\Lambda_aQ_a(z)}B_a(z)\,dv_\alpha(z)\right|
\le C_r\Lambda_a^{-1}
=C_r(1-|a|^2)^\gamma.
\]
Together with \eqref{eq:osc-average-ratio}, this proves
\eqref{eq:osc-average-rate}.  In particular, the main theorem already
implies that $T_{\phi_\gamma}$ is compact.

\smallskip
\noindent\emph{Step 2: radial diagonalization.}
We now exploit the radial form of \eqref{eq:osc-symbol}.  Write
\[
A^2_\alpha(\B_n)=\bigoplus_{\ell=0}^\infty \mathcal H_\ell,
\]
where $\mathcal H_\ell$ denotes the finite-dimensional space of
homogeneous holomorphic polynomials of degree $\ell$.  For
$f,g\in\mathcal H_\ell$, polar coordinates show that multiplication by a
radial symbol changes only the radial integral and leaves the spherical
inner product unchanged.  If $f\in\mathcal H_\ell$ and
$g\in\mathcal H_j$ with $j\ne\ell$, invariance under
$z\mapsto e^{i\theta}z$ gives orthogonality.  Hence
\begin{equation}\label{eq:radial-diagonal}
T_{\phi_\gamma}|_{\mathcal H_\ell}=\lambda_\ell I_{\mathcal H_\ell},
\end{equation}
where
\begin{equation*}
\lambda_\ell
=\frac{\displaystyle
\int_0^1 e^{i(1-r^2)^{-\gamma}}
 r^{2\ell+2n-1}(1-r^2)^\alpha\,dr}
{\displaystyle
\int_0^1 r^{2\ell+2n-1}(1-r^2)^\alpha\,dr}.
\end{equation*}
After the change of variables $t=1-r^2$, the factor $1/2$ cancels and
\begin{equation*}
\lambda_\ell
=\frac{I_\ell}{B(\alpha+1,\ell+n)},
\qquad
I_\ell=\int_0^1 e^{it^{-\gamma}}
 t^\alpha(1-t)^{\ell+n-1}\,dt.
\end{equation*}
The beta-function denominator satisfies
\begin{equation}\label{eq:beta-denominator}
B(\alpha+1,\ell+n)
\asymp (1+\ell)^{-\alpha-1}.
\end{equation}

\smallskip
\noindent\emph{Step 3: super-polynomial decay of the radial eigenvalues.}
Introduce
\[
\mathcal D_\gamma f(t)
=\frac{d}{dt}\bigl(t^{\gamma+1}f(t)\bigr).
\]
Since
\[
\frac{d}{dt}e^{it^{-\gamma}}
=-i\gamma t^{-\gamma-1}e^{it^{-\gamma}},
\]
we have
\[
e^{it^{-\gamma}}
=\frac{i}{\gamma}t^{\gamma+1}
\frac{d}{dt}e^{it^{-\gamma}}.
\]
For a fixed integer $N\ge1$ and all sufficiently large $\ell$, repeated
integration by parts is legitimate and all boundary terms vanish.  At
$t=0$ this follows from $\alpha>-1$ and $\gamma>0$; after $j$ steps every
boundary amplitude has at least the factor
$t^{\alpha+(j+1)\gamma+1}$.  At $t=1$ it follows from the zero of
$(1-t)^{\ell+n-1}$, whose order tends to infinity with $\ell$.  Thus
\begin{equation}\label{eq:IBP-N}
I_\ell
=\left(-\frac{i}{\gamma}\right)^N
\int_0^1 e^{it^{-\gamma}}
\mathcal D_\gamma^N
\bigl[t^\alpha(1-t)^{\ell+n-1}\bigr]dt.
\end{equation}

We next record the elementary structure of the differentiated amplitude.
An induction based on
\[
\mathcal D_\gamma\bigl[t^a(1-t)^b\bigr]
=(a+\gamma+1)t^{a+\gamma}(1-t)^b
-bt^{a+\gamma+1}(1-t)^{b-1}
\]
gives
\begin{equation}\label{eq:Dgamma-expansion}
\mathcal D_\gamma^N
\bigl[t^\alpha(1-t)^M\bigr]
=\sum_{j=0}^N c_{N,j}(M)
 t^{\alpha+N\gamma+j}(1-t)^{M-j},
\end{equation}
where $c_{N,j}(M)$ is a polynomial in $M$ of degree at most $j$.
Consequently
\begin{equation*}
|c_{N,j}(M)|\le C_N(1+M)^j,
\qquad 0\le j\le N.
\end{equation*}
Taking $M=\ell+n-1$ in \eqref{eq:Dgamma-expansion} and using
\eqref{eq:IBP-N}, we find
\begin{align}
|I_\ell|
&\le C_N\sum_{j=0}^N(1+\ell)^j
 B(\alpha+N\gamma+j+1,\ell+n-j)\notag\\
&\le C_N(1+\ell)^{-\alpha-1-N\gamma}.
\label{eq:Iell-decay}
\end{align}
Here the second line follows from the standard gamma-ratio estimate
$B(A,\ell+B)\asymp \ell^{-A}$ for fixed $A>0$ and fixed $B$.
Combining \eqref{eq:Iell-decay} with \eqref{eq:beta-denominator} yields
\eqref{eq:eigen-superpoly}.  Enlarging $C_N$ absorbs the finitely many
small values of $\ell$ for which the repeated integration by parts was
not used.
\smallskip
\noindent\emph{Step 4: Schatten membership and the Carleson contrast.}
The multiplicity of the eigenvalue $\lambda_\ell$ is
\begin{equation}\label{eq:multiplicity}
d_\ell=\dim\mathcal H_\ell
=\binom{\ell+n-1}{n-1}
\asymp(1+\ell)^{n-1}.
\end{equation}
Since \eqref{eq:radial-diagonal} makes $T_{\phi_\gamma}$ a normal
diagonal operator, for every $p>0$,
\[
T_{\phi_\gamma}\in\mathcal S_p
\quad\Longleftrightarrow\quad
\sum_{\ell=0}^\infty d_\ell|\lambda_\ell|^p<\infty.
\]
Given $p>0$, choose $N$ so large that $N\gamma p>n$.  By
\eqref{eq:eigen-superpoly} and \eqref{eq:multiplicity},
\[
\sum_{\ell=0}^\infty d_\ell|\lambda_\ell|^p
\le C_N\sum_{\ell=0}^\infty
(1+\ell)^{n-1-N\gamma p}<\infty.
\]
This proves \eqref{eq:all-schatten}, including the quasi-Schatten range
$0<p<1$.

For completeness, let $s_j(T_{\phi_\gamma})$ denote the singular values
in decreasing order. Truncation to homogeneous degrees at most $L$ has rank
$\sum_{\ell=0}^Ld_\ell$ and operator-norm error at most
$\sup_{\ell>L}|\lambda_\ell|$. Since
\[
\sum_{\ell=0}^L d_\ell=\binom{L+n}{n}\asymp L^n,
\]
\eqref{eq:eigen-superpoly} implies that for every $M>0$ there exists
$C_M$ such that
\[
s_j(T_{\phi_\gamma})\le C_M j^{-M},
\qquad j\ge1.
\]
Thus the singular values decay faster than every power.

Finally, $|\phi_\gamma|\equiv1$, and hence for every $a$ and $r$,
\[
\frac{(|\phi_\gamma|dv_\alpha)(D(a,r))}
{v_\alpha(D(a,r))}=1.
\]
Thus $|\phi_\gamma|\,dv_\alpha=dv_\alpha$ is not a vanishing
Carleson measure.
\end{proof}

\begin{remark}\label{rem:carleson-distinction}
For the symbols in Proposition~\ref{prop:oscillatory-symbol}, the
average of $|\phi_\gamma|$ equals one on every ball, whereas the
average of $\phi_\gamma$ tends to zero and the Toeplitz singular
values decay faster than every power. Thus a condition on the
modulus cannot characterize compactness for complex-valued symbols.
\end{remark}

\section{Further examples}
\label{lim:section}

The following examples examine the radius, the boundary approach,
and the metric in the averaging criterion.

\subsection{A radius depending on the center}
\label{lim:variable-radius}

A sequence of fixed-radius boundary tests is different from a single
test whose radius depends on the center. In particular, the fixed
radius in Theorem~\ref{thm:psc-fixed-radius} cannot be replaced by an
arbitrary positive function tending to zero at the boundary.

\begin{proposition}
\label{lim:shrinking-counterexample}
On the unit ball $\mathbb B_n$, set $\delta_{\mathbb B}(z)=1-|z|^2$ and fix
$\gamma>0$. The symbol
\[
 \phi_\gamma(z)=\exp\bigl(i\delta_{\mathbb B}(z)^{-\gamma}\bigr)
\]
has $T_{\phi_\gamma}\in\mathcal S_p$ for every $p>0$, whereas,
for $R(a)=\delta_{\mathbb B}(a)^{2\gamma}$,
\[
 \mathcal A_{R(a)}^{\mathbb B_n}\phi_\gamma(a)
           -\phi_\gamma(a)\longrightarrow0,
 \qquad
 \left|\mathcal A_{R(a)}^{\mathbb B_n}\phi_\gamma(a)\right|
           \longrightarrow1
 \quad(a\to\partial\mathbb B_n).
\]
The average uses the same fixed normalization of the Bergman metric
as Theorem~\ref{thm:psc-fixed-radius}.
\end{proposition}

\begin{proof}
The Schatten assertion is Proposition~\ref{prop:oscillatory-symbol}.
For the assertion about the mean, the ball Bergman metric satisfies
\[
 |d\log\delta_{\mathbb B}|_{g_{\mathbb B_n}}\le C,
 \qquad
 |d(\delta_{\mathbb B}^{-\gamma})|_{g_{\mathbb B_n}}
                 \le C_\gamma\delta_{\mathbb B}^{-\gamma}.
\]
On a metric ball of radius $R(a)\le1$, the first estimate makes
$\delta_{\mathbb B}$ comparable to $\delta_{\mathbb B}(a)$, with constants independent of
$a$. Integrating the second estimate along a minimizing metric
geodesic therefore gives
\[
 \sup_{w\in D_{R(a)}(a)}
    |\delta_{\mathbb B}(w)^{-\gamma}-\delta_{\mathbb B}(a)^{-\gamma}|
 \le C_\gamma R(a)\delta_{\mathbb B}(a)^{-\gamma}
 =C_\gamma\delta_{\mathbb B}(a)^\gamma\longrightarrow0.
\]
Since $|e^{it}-e^{is}|\le|t-s|$, the same bound controls
$\sup_{D_{R(a)}(a)}|\phi_\gamma-\phi_\gamma(a)|$.
It controls the difference between the normalized Euclidean mean
and $\phi_\gamma(a)$ as well. The last assertion follows from
$|\phi_\gamma(a)|=1$.
\end{proof}

This example concerns necessity: even an operator lying in every
Schatten class can fail the vanishing condition at a center-dependent
shrinking radius. There is no conflict with
Corollary~\ref{pd:shrinking-radii}, which first takes a boundary
limit at each fixed radius $R_k$, separately for every $k$.
Its quantifiers are
\[
 \text{for every }k,\qquad
 \lim_{a\to\partial\mathbb D^n}A_{R_k}f(a)=0,
\]
not a simultaneous limit with $k$ depending on $a$.

\subsection{The entire topological boundary is required}
\label{lim:whole-boundary}

In the polydisc results, $a\to\partial\mathbb D^n$ means
$\max_j|a_j|\to1$. It does not mean that every coordinate must
approach its unit circle. Testing only the latter approach would
miss noncompact operators.

\begin{proposition}
\label{lim:distinguished-boundary}
Let $n\ge2$ and
\[
 f(z)=1-|z_1|^2\qquad(z\in\mathbb D^n).
\]
Then $T_f$ is not compact. Nevertheless, for every fixed positive
geodesic radius $R$ and every fixed positive box radius vector
$\mathbf r$,
\[
 A_R f(a)\longrightarrow0,
 \qquad A_{\mathbf r}^{\mathrm{box}}f(a)\longrightarrow0
 \quad\text{whenever }|a_j|\to1\text{ for all }j.
\]
\end{proposition}

\begin{proof}
For a fixed invariant radius, every point in the geodesic ball, or
in the product box, satisfies
\[
 1-|z_1|^2\le C(1-|a_1|^2),
\]
where $C$ depends only on the chosen radii and metric normalization.
For geodesic balls this follows by projecting the product metric
onto the first disc; for boxes it is the one-dimensional
fixed-radius height comparison. Nonnegativity of $f$ gives the
asserted limits.

On the other hand, the tensor product identification of the
Bergman space gives
\[
 T_f=T_{1-|z|^2}\otimes I
 \quad\text{on }A^2(\mathbb D)\otimes A^2(\mathbb D^{n-1}).
\]
The first factor is nonzero, while the second space is infinite
dimensional. If $v$ is a unit vector with
$T_{1-|z|^2}v\ne0$ and $(e_k)$ is an orthonormal sequence in the
second space, the vectors $T_f(v\otimes e_k)$ are mutually
orthogonal and have the same positive norm. Hence $T_f$ is not
compact.
\end{proof}

Thus the whole-boundary condition in Theorem~\ref{pd:box} cannot
be replaced by simultaneous coordinatewise approach to the
distinguished boundary. Independently,
Theorem~\ref{pd:counterexample} shows that fixed product boxes and
fixed Bergman geodesic balls have different detecting properties,
even when both boundary limits are taken along the entire
boundary.

\subsection{Fixed Euclidean balls do not give the same criterion}
\label{lim:euclidean-balls}

Although the averages use Euclidean volume, their sets are Bergman
metric balls. Averages over Euclidean balls of a fixed radius fail
necessity even for compactly supported symbols.

\begin{proposition}
\label{lim:euclidean-counterexample}
Let $\Omega\Subset\mathbb C^n$ be a domain, let
$p\in\partial\Omega$, and fix $r>0$. There exists a nonzero
$f\in C_c^\infty(\Omega)$, $f\ge0$, for which $T_f$ is compact,
but
\[
 \begin{gathered}
E_r f(a):=\frac{1}{V(\Omega\cap B_{\mathrm{Euc}}(a,r))}
      \int_{\Omega\cap B_{\mathrm{Euc}}(a,r)}f\,dV
 \\
 \longrightarrow c>0\qquad(a\to p,\ a\in\Omega).
 \end{gathered}
\]
\end{proposition}

\begin{proof}
Choose a nonzero nonnegative smooth bump with
$\operatorname{supp}f\Subset\Omega\cap B_{\mathrm{Euc}}(p,r/2)$.
For all $a$ sufficiently close to $p$, the integration ball contains
the support of $f$, so its numerator is the positive constant
$\int_\Omega f\,dV$.
Euclidean spheres have Lebesgue measure zero, and dominated
convergence gives
\[
 V(\Omega\cap B_{\mathrm{Euc}}(a,r))
    \longrightarrow V(\Omega\cap B_{\mathrm{Euc}}(p,r))>0.
\]
This proves the asserted positive limit. A bounded compactly
supported symbol defines a compact Toeplitz operator: restriction
of bounded subsets of $A^2(\Omega)$ to its compact support is
compact in $L^2$, by local estimates and Montel's theorem.
\end{proof}

\subsection{What the domain hypotheses establish}
\label{lim:domain-hypotheses}

The scaling argument requires continuity of the second boundary
derivatives, without a H\"older modulus; see
\eqref{eq:geom-C2-estimate}. Appendix~\ref{app:boundary-regularity}
gives a strongly convex $C^2$ domain whose boundary is not
$C^{2,\alpha}$ for any $\alpha>0$.

The optimal domain class remains open. The polydisc counterexample
rules out an extension to all bounded pseudoconvex domains, but does
not settle the criterion on smooth weakly pseudoconvex domains.

\appendix
\section{A Fourier--Hermite proof of the model rigidity}
\label{app:fourier-rigidity}

We give an alternative proof of
Proposition~\ref{sp:eigenfunction-rigidity} using the full
special-Hermite expansion and a Banach-valued Euler equation.

\subsection{A Banach-valued Euler equation}

\begin{lemma}
\label{sp:banach-euler}
Let $E$ be a Banach space, $A,B\in\mathcal L(E)$, $\nu>n^2$, and
$m>n/2$. Suppose $f\in C^2((0,\rho_0);E)$ solves
\begin{equation}
 \rho^2f''+(1-n)\rho f'+\frac\nu4 f
       =\rho Af+\rho^2Bf
 \label{sp:banach-equation}
\end{equation}
and $\|f(\rho)\|_E\le C\rho^m$ for $0<\rho<\rho_0$.
Then $f=0$ on $(0,\rho_0)$.
\end{lemma}

\begin{proof}
Set
\[
 \tau=\frac12\sqrt{\nu-n^2}>0,
 \qquad s_\pm=\frac n2\pm i\tau.
\]
These are the roots of $s(s-n)+\nu/4=0$.
The Green kernel for
$\partial_\rho^2+(1-n)\rho^{-1}\partial_\rho+\nu/(4\rho^2)$
with no homogeneous component at zero is
\begin{equation*}
 G(\rho,t)=
 \frac{\rho^{s_+}t^{s_-}-\rho^{s_-}t^{s_+}}
      {(s_+-s_-)t^{n-1}},\qquad 0<t<\rho.
\end{equation*}
Indeed, $G(t,t)=0$ and $\partial_\rho G(t,t)=1$.
Writing $t=\rho a$ gives $G(\rho,\rho a)=\rho g(a)$, where
\[
 |g(a)|\le C_\nu a^{1-n/2},\qquad 0<a<1.
\]
Consequently
\begin{equation}
 \int_0^1|g(a)|a^{m-1}\,da+
 \int_0^1|g(a)|a^m\,da<\infty.
 \label{sp:green-integrability}
\end{equation}

The Bochner integral
\[
 I(\rho)=\int_0^\rho G(\rho,t)
                   \bigl(t^{-1}Af(t)+Bf(t)\bigr)\,dt
\]
is well defined, solves the inhomogeneous normalized equation, and is
$O(\rho^{m+1})$ as $\rho\downarrow0$.
Thus $f-I$ solves the homogeneous Euler equation and is $O(\rho^m)$.
Pairing with any element of $E^*$ reduces this assertion to a scalar
linear combination of $\rho^{s_+}$ and $\rho^{s_-}$.
Since $m>\operatorname{Re}s_\pm=n/2$, that combination must vanish.
The Hahn--Banach theorem therefore gives
\begin{equation}
 f(\rho)=\int_0^\rho G(\rho,t)
                   \bigl(t^{-1}Af(t)+Bf(t)\bigr)\,dt.
 \label{sp:volterra}
\end{equation}
For $0<r<\rho_0$, let
$N_r=\sup_{0<\rho\le r}\rho^{-m}\|f(\rho)\|_E$.
Equations \eqref{sp:green-integrability}--\eqref{sp:volterra} imply
\[
 N_r\le \bigl(C_1\|A\|r+C_2\|B\|r^2\bigr)N_r.
\]
For small $r$, this forces $N_r=0$.
Uniqueness for the Banach-valued ordinary differential equation on
closed subintervals of $(0,\rho_0)$ then gives the result.
\end{proof}

\subsection{Alternative proof of the eigenfunction consequence}

\begin{proof}[Alternative proof of Proposition~\ref{sp:eigenfunction-rigidity}]
Use an isometry sending $\zeta$ to infinity in the Siegel domain
\[
 \mathcal U^n=\{(z,w)\in\mathbb C^{n-1}\times\mathbb C:
                    \rho=\operatorname{Im}w-|z|^2>0\}.
\]
The height is a positive constant multiple of $\rho$; absorbing this
constant into $C,p_0$, the assumption becomes
\begin{equation}
 |u(z,t,\rho)|\le C\rho^m
 \quad\bigl((z,t)\in\mathbb C^{n-1}\times\mathbb R,
                  \ 0<\rho\le\rho_0\bigr),
 \label{sp:siegel-bound}
\end{equation}
where $t=\operatorname{Re}w$.
Put
\[
 X_j=\partial_{x_j}+2y_j\partial_t,
 \quad Y_j=\partial_{y_j}-2x_j\partial_t,
 \quad T=\partial_t,
 \quad L_0=-\frac14\sum_{j=1}^{n-1}(X_j^2+Y_j^2).
\]
In the chosen normalization,
\begin{equation*}
 dV=c_n\rho^{-n-1}\,dz\,dt\,d\rho,
 \qquad
 \Delta_X=4\rho\bigl[\rho(\partial_{\rho\rho}+T^2)
                     -L_0-(n-1)\partial_\rho\bigr].
\end{equation*}
The eigenfunction equation is therefore
\begin{equation}
 \rho^2u_{\rho\rho}-(n-1)\rho u_\rho+\frac\nu4u
                  =\rho L_0u-\rho^2T^2u.
 \label{sp:siegel-equation}
\end{equation}
For every closed height interval contained in $(0,\rho_0)$,
interior elliptic estimates, followed by Heisenberg translations, give
uniform bounds on $\partial_\rho^j u$ for each fixed $j$.
In particular, all Fourier coefficients below depend smoothly on $\rho$
in the stated negative Sobolev spaces, not only as distributions.

\emph{Nonzero central frequencies.}
First suppose $n\ge2$ and Fourier transform in $t$, with dual variable
$\lambda$. At $\lambda\ne0$, the transformed operator $L_0^\lambda$
is the twisted Laplacian. It has a complete special-Hermite basis
$\theta_j^\lambda$ of Schwartz functions with eigenvalues
\[
 L_0^\lambda\theta_j^\lambda
          =c_j|\lambda|\theta_j^\lambda,\qquad c_j>0.
\]
Here the index enumerates the full basis, including multiplicities.
The basis can be chosen smoothly for $\lambda>0$ and for $\lambda<0$.
On compact subsets of $\mathbb R\setminus\{0\}$, its Schwartz seminorms
and every fixed number of its $\lambda$ derivatives are uniformly bounded,
for each fixed $j$. Its expansions converge in the Schwartz topology,
locally uniformly with respect to $\lambda$ away from zero.
See Folland~\cite{Folland} and Thangavelu~\cite{Thangavelu}.
These standard parameter-dependent properties also follow by dilation
from the special-Hermite basis at $\lambda=1$ and at $\lambda=-1$.

Fix $\chi\in C_c^\infty(\mathbb R\setminus\{0\})$, put
$K=\operatorname{supp}\chi$, and let $f_j(\rho)$ be the distribution in
$\lambda$ obtained by pairing $\chi\widehat u$ in the horizontal variable
with $\theta_j^\lambda$.
More precisely, using bilinear distribution duality, set
\[
 \langle f_j(\rho),\psi\rangle
   =\bigl\langle\widehat u(\cdot,\cdot,\rho),
         \chi(\lambda)\psi(\lambda)
                  \overline{\theta_j^\lambda(z)}\bigr\rangle.
\]
The transpose of $L_0^\lambda$ is $L_0^{-\lambda}$, and
$L_0^{-\lambda}\overline{\theta_j^\lambda}
      =c_j|\lambda|\overline{\theta_j^\lambda}$.
Thus the eigenvalue reduction is valid in this distribution pairing.
For completeness, the uniform bound required for this distributional
pairing is as follows. If $\psi\in C_c^\infty(\mathbb R)$, the inverse
Fourier kernel of
$\chi(\lambda)\psi(\lambda)\overline{\theta_j^\lambda(z)}$
has $L^1(dz\,dt)$ norm at most
\begin{equation*}
 C_{j,\chi}\|\psi\|_{H^2(\mathbb R)}.
\end{equation*}
For $|t|\le1$ this follows directly from the uniform Schwartz bounds;
for $|t|>1$, integrate twice in $\lambda$ and use $t^{-2}$ decay.
Equation \eqref{sp:siegel-bound} consequently gives
\begin{equation}
 \|f_j(\rho)\|_{H^{-2}(\mathbb R)}
            \le C_{j,\chi}\rho^m.
 \label{sp:coefficient-sobolev-bound}
\end{equation}
The same test-kernel estimate, applied to the height derivatives on
closed height intervals, proves
$f_j\in C^2((0,\rho_0);H^{-2}(\mathbb R))$.

Take as $E$ the closed subspace of $H^{-2}(\mathbb R)$ consisting of
distributions supported in $K$.
Multiplication on $E$ by $c_j|\lambda|$ and by $\lambda^2$ is bounded:
extend these functions smoothly from a neighborhood of $K$ to compactly
supported smooth functions on $\mathbb R$.
The coefficient equation from \eqref{sp:siegel-equation} is
\[
 \rho^2f_j''+(1-n)\rho f_j'+\frac\nu4 f_j
       =\rho\,c_j|\lambda|f_j+\rho^2\lambda^2f_j.
\]
Lemma~\ref{sp:banach-euler} and
\eqref{sp:coefficient-sobolev-bound} give $f_j=0$.
Completeness in the Schwartz topology, followed by duality, implies
$\chi\widehat u=0$.
Since $\chi$ is arbitrary, the partial Fourier transform of $u$ in $t$
is supported at $\lambda=0$ for every $0<\rho<\rho_0$.

For $n=1$ there is no horizontal variable and $L_0=0$.
The same argument is applied directly to $\chi\widehat u$,
with $A=0$ and $B$ multiplication by $\lambda^2$; no special-Hermite
expansion is used in this case.

The point-support theorem for distributions now implies $Tu=0$ at these
heights. One can apply it after pairing with any compactly supported
horizontal test function: the resulting function of $t$ is bounded,
has Fourier transform supported at zero, and is therefore a constant.
This avoids any assertion about multiplication of singular coefficients
by point-supported distributions.

\emph{Nonzero horizontal frequencies.}
For $n\ge2$, $Tu=0$ reduces $L_0$ to $-\Delta_z/4$, so that
\eqref{sp:siegel-equation} becomes
\[
 \rho^2u_{\rho\rho}-(n-1)\rho u_\rho+\frac\nu4u
                        =-\frac\rho4\Delta_z u.
\]
Fourier transform in $z\in\mathbb R^{2n-2}$ and localize by a compactly
supported cutoff away from $\xi=0$.
The localized transform belongs to $H^{-N}(\mathbb R^{2n-2})$, with norm
$O(\rho^m)$, for any sufficiently large fixed integer $N$; this follows
by bounding the $L^1$ norm of the inverse Fourier test kernel by its
$H^N$ norm. It is $C^2$ in height by the same strip estimates.
Its equation has the form \eqref{sp:banach-equation} with
$A$ multiplication by $|\xi|^2/4$ and $B=0$, on the closed subspace
of distributions supported in the localization set.
Lemma~\ref{sp:banach-euler} makes every such localization zero.
The horizontal Fourier transform is thus supported at the origin.
The point-support theorem and boundedness in $z$ imply that $u$ is also
independent of $z$.

We have proved, for every $n\ge1$, that $u=u(\rho)$ when
$0<\rho<\rho_0$.
It solves
\[
 \rho^2u''+(1-n)\rho u'+\frac\nu4u=0.
\]
Its two indicial powers have real part $n/2$, whereas
$u=O(\rho^m)$ with $m>n/2$.
Thus $u=0$ on this open height region.
Unique continuation for the elliptic equation
$(\Delta_X+\nu)u=0$ gives $u=0$ on the connected manifold $X$.
\end{proof}

\section{An auxiliary boundary-regularity example}
\label{app:boundary-regularity}

The following example has continuous second boundary derivatives
but no H\"older modulus of positive order.

\begin{proposition}[A strongly convex $C^2$ domain with no H\"older second derivative]
\label{lim:bare-C2}
Let $n\ge2$ and $\varepsilon>0$. Define
\[
 \begin{gathered}
\omega(0)=0,\qquad
 \omega(s)=\frac1{\log(e+s^{-2})}\quad(s>0),\\
 h(x)=\int_0^{|x|}(|x|-s)\omega(s)\,ds.
 \end{gathered}
\]
The domain
\[
 \Omega_\varepsilon
 =\{z\in\mathbb C^n:|z|^2+\varepsilon h(\operatorname{Re}z_1)<1\}
\]
has $C^2$ strongly convex, hence strongly pseudoconvex, boundary.
Its boundary is not $C^{2,\alpha}$ for any $0<\alpha\le1$.
Nevertheless all the equivalences of
Theorem~\ref{thm:psc-fixed-radius} hold on this domain.
\end{proposition}

\begin{proof}
The function $h$ is even, nonnegative, and $C^2$ on $\mathbb R$,
with $h(0)=h'(0)=0$ and $h''(x)=\omega(|x|)\ge0$.
For $r(z)=|z|^2+\varepsilon h(\operatorname{Re}z_1)-1$, one has
\[
 (r_{j\bar k})=\operatorname{diag}
 \left(1+\frac\varepsilon4\omega(|\operatorname{Re}z_1|),
                  1,\ldots,1\right)\ge I.
\]
The real Hessian is at least $2I$. On $r=0$, writing
$x=\operatorname{Re}z_1$, the radial derivative is
\[
 dr(z)[z]=2|z|^2+\varepsilon xh'(x)>0.
\]
The domain is bounded because $h\ge0$, and its boundary is therefore
$C^2$ and strongly convex.

It remains to check the regularity of the boundary itself, rather
than just that of this particular defining function. At
$p=(0,1,0,\ldots,0)$, the boundary is a graph in the
$\operatorname{Re}z_2$ direction. Its restriction obtained by putting
$z_1=x\in\mathbb R$ and all other tangential coordinates equal to
zero is
\[
 u(x)=\sqrt{1-q(x)},\qquad q(x)=x^2+\varepsilon h(x).
\]
In particular $u''(0)=-1$, and for small $x\ne0$,
\begin{align*}
 -u''(x)-1
 &=\left(\frac1{\sqrt{1-q(x)}}-1\right)
   +\frac{\varepsilon\omega(|x|)}{2\sqrt{1-q(x)}}
   +\frac{q'(x)^2}{4(1-q(x))^{3/2}}\\
 &\ge\frac\varepsilon2\omega(|x|).
\end{align*}
Since $\omega(|x|)/|x|^\alpha\to\infty$ for every $\alpha>0$,
this graph is not $C^{2,\alpha}$ for any $0<\alpha\le1$.
The same is therefore true of the boundary. The final assertion is
the application of Theorem~\ref{thm:psc-fixed-radius}.
\end{proof}

\section*{Acknowledgments}

\noindent\textbf{Funding.} The authors were partially supported by the
National Natural Science Foundation of China (Grant numbers
12171484 and 12571140), the Natural Science Foundation of Hunan Province
(Grant number 2023JJ20056), the Science and Technology Innovation Program
of Hunan Province (Grant number 2023RC3028), the Central South University
Innovation-Driven Research Programme (Grant number 2023CXQD032), and
the Hunan Basic Science Research Center for Mathematical Analysis
(2024JC2002).

\smallskip
\noindent\textbf{AI Statement.}
The author used OpenAI’s ChatGPT (GPT-5.6 Sol) for language editing,
LaTeXformatting,
and limited assistance with local mathematical reasoning. The proof strategy
and mathematical
development are the author’ own; all mathematical content was
independently written and checked by the author, who take full responsibility for
it.

\end{document}